\documentclass{article}

\usepackage{microtype}
\usepackage{graphicx}
\usepackage{subfigure}
\usepackage{booktabs} % for professional tables

\usepackage{hyperref}

 \usepackage[accepted]{icml2024}

\usepackage{amsmath}
\usepackage{amssymb}
\usepackage{mathtools}
\usepackage{amsthm}

\usepackage[capitalize,noabbrev]{cleveref}

\theoremstyle{plain}
\newtheorem{theorem}{Theorem}
\newtheorem{proposition}{Proposition}
\newtheorem{lemma}{Lemma}

\theoremstyle{definition}
\newtheorem{definition}{Definition}
\newtheorem{assumption}{Assumption}
\theoremstyle{remark}
\newtheorem{remark}{Remark}

\usepackage[textsize=tiny]{todonotes}

\icmltitlerunning{Online Optimization Perspective on First-Order and Zero-Order Decentralized Nonsmooth Nonconvex Stochastic Optimization}

\begin{document}

\twocolumn[
\icmltitle{An Online Optimization Perspective on First-Order and Zero-Order Decentralized Nonsmooth Nonconvex Stochastic Optimization}



\icmlsetsymbol{equal}{*}

\begin{icmlauthorlist}
\icmlauthor{Emre Sahinoglu}{yyy}
\icmlauthor{Shahin Shahrampour}{yyy}
\end{icmlauthorlist}

\icmlaffiliation{yyy}{Department of Mechanical and Industrial Engineering, Northeastern University, Boston, MA 02115}

\icmlcorrespondingauthor{Emre Sahinoglu}{sahinoglu.m@northeastern.edu}
\icmlcorrespondingauthor{Shahin Shahrampour}{s.shahrampour@northeastern.edu}

\icmlkeywords{Optimization, Decentralized Optimization}

\vskip 0.3in
]



\printAffiliationsAndNotice{}  % leave blank if no need to mention equal contribution

\newcommand{\goldstat}{($\delta,\epsilon$)-stationary }

\begin{abstract}
We investigate the finite-time analysis of finding \goldstat points for nonsmooth nonconvex objectives in decentralized stochastic optimization. A set of agents aim at minimizing a global function using only their local information by interacting over a network. We present a novel algorithm, called Multi Epoch Decentralized Online Learning (ME-DOL), for which we establish the sample complexity in various settings. First, using a recently proposed online-to-nonconvex technique, we show that our algorithm recovers the optimal convergence rate of smooth nonconvex objectives. We then extend our analysis to the nonsmooth setting, building on properties of randomized smoothing and Goldstein-subdifferential sets. We establish the sample complexity of $O(\delta^{-1}\epsilon^{-3})$, which to the best of our knowledge is the first finite-time guarantee for decentralized nonsmooth nonconvex stochastic optimization in the first-order setting (without weak-convexity), matching its optimal centralized counterpart. We further prove the same rate for the zero-order oracle setting without using variance reduction.
\end{abstract}

\newcommand{\norm}[1]{\left \lVert #1 \right \rVert}
\newcommand{\gradnormdelta}[2]{\norm {\nabla #1 (#2) }_\delta }

\newcommand{\complexityns}[2]{$O(\epsilon^{#1}\delta^{#2})$ }
\newcommand{\oracle}[1]{\mathcal{O}_{#1}}
\newcommand{\un}{\texttt{unif}}
\newcommand{\bl}{\color{blue}}
\newcommand{\rd}{\color{red}}

\newcommand{\innerproduct}[1]{\langle #1 \rangle}

\newcommand{\finitesum}[2]{\sum_{#1=1}^#2}

\newcommand{\lhslocal}{\finitesum{t}{T} \finitesum{i}{n} \nabla f^i(w_{t,i}^k)}

\newcommand{\lhsglobal}{\finitesum{t}{T} \finitesum{i}{n} \nabla f(w_{t,i}^k)}

\newcommand{\lhsglobalf}[1]{\finitesum{t}{T} \finitesum{i}{n} \nabla #1(w_{t,i}^k)}

\newcommand{\probexpbig}[1]{\mathbb{E}\left[ #1 \right]}

\newcommand{\probexp}[1]{\mathbb{E}[ #1 ]}

\newcommand{\epsterm}{\probexpbig{ \frac{1}{K} \sum_{k=1}^K\norm{\frac{1}{nT} \lhslocal}}}

\newcommand{\epstermglobal}{\probexpbig{ \frac{1}{K} \sum_{k=1}^K\norm{\frac{1}{nT} \lhsglobal}}}

\newcommand{\epstermglobalf}[1]{\probexpbig{ \frac{1}{K} \sum_{k=1}^K\norm{\frac{1}{nT} \lhsglobalf{#1}}}}

\newcommand{\regretterm}{\frac{\probexp{\sum_{t=1}^{T}\innerproduct{ \bar{g}_{t}^k , \bar{\Delta}_t^k - u^k }}}{DT}}

\newcommand{\discrepancyterm}{\frac{\mathbb{E}[ \sum_{t=1}^{T}\innerproduct{\bar{\Delta}_t^k ,\tilde{\nabla}_t^k - \bar{\nabla}_t^k }]}{DT}}

\newcommand{\suboptterm}{\frac{f(\bar{x}_0) - \mathbb{E}[f(\bar{x}_{TK})]}{DTK}}

\newcommand{\zeroordersecbound}{16\sqrt{2\pi} d L^2}

\newcommand{\geoconstant}{c}

\section{Introduction}
At the heart of many practical machine learning problems, we must deal with nonconvex optimization of nonsmooth objective functions. Examples include training neural networks with ReLU activation functions, blind deconvolution, sparse dictionary learning, and robust phase retrieval. Despite the significant practical success of such schemes, the vast majority of prior work in theoretical analysis of nonsmooth nonconvex optimization focused on asymptotic convergence results \cite{rockafellar2009variational,clarke2008nonsmooth}. More recently, {\it finite-time} analysis of this class of problems has attracted significant attention \cite{jordan2023deterministic,Majewski2018AnalysisON,davis2019stochastic,daniilidis2020pathological,tian2022finitetime}. 

On the other hand, decentralization is a crucial mechanism to scale up optimization problems. For nonsmooth objectives, though the class of convex problems is well-understood in decentralized optimization \cite{nedic09distributed,scaman2018optimal}, the characterization of optimal finite-time rates for nonconvex problems has remained elusive (except for weakly-convex problems \cite{Chen_Distributed_weaklyconvex}). In the present work, we address the {\it finite-time} analysis of {\it decentralized nonsmooth nonconvex stochastic} optimization. 

We consider a decentralized optimization problem where a group of $n$ agents aim at minimizing a global function. However, each agent has limited information about this global objective and interacts with its neighbors to solve the global problem, formulated in the following form
\begin{align}\label{eq:global}
  \displaystyle \min_{x \in \mathbb{R}^d} \Big\{ f(x) = \frac{1}{n} \sum_{i=1}^{n} f^{i}(x) \Big\}.  
\end{align}
Local functions $f^i$ are in the form of $f^{i}(x) = \mathbb{E}_{\xi_i \sim \mathcal{D}_i} [F^i (x,\xi_i)]$, where $F^i(x,\xi_i)$ are stochastic with random index $\xi_i$, and $\xi_i$ corresponds to a data sample from local dataset of agent $i$. We assume that the local functions are nonconvex and Lipschitz continuous but do {\it not} necessarily have Lipschitz continuous gradients, i.e., they are nonsmooth.

In an optimization problem, a tractable optimality criterion is required for finite-time convergence guarantees. For nonsmooth nonconvex objectives, $\epsilon$-stationarity cannot be guaranteed in finite time \cite{kornowski2022oracle,zhang2020complexity}. Instead, the notion of ($\delta,\epsilon$)-stationarity  is a tractable criterion \cite{zhang2020complexity}, where we seek vectors with norm less than $\epsilon$ among the convex hull of the subdifferential set of a ball with radius $\delta$ (see Definition \ref{def:goldstein} for exact mathematical definition). The goal of this paper is to identify a \goldstat point of the global function $f$ when agents have access to either the first-order oracle (i.e., $\nabla F^i(\cdot,\xi_i)$) or the zero-order oracle (i.e., $F^i(\cdot,\xi_i)$).

\subsection{Contributions}
In this paper, we address the finite-time analysis of decentralized nonsmooth nonconvex stochastic optimization. We present a novel algorithm, called \textbf{M}ulti \textbf{E}poch \textbf{D}ecentralized \textbf{O}nline \textbf{L}earning (ME-DOL), for which we establish the sample complexity in various settings. Our contributions are three-fold.

\begin{itemize}
\item We adopt the online-to-nonconvex conversion technique of \citet{cutkosky2023optimal} to streamline the finite-time analysis of ME-DOL for decentralized nonsmooth nonconvex optimization. First, for smooth objectives, we prove that the complexity of finding \goldstat points is $O(\delta^{-1}\epsilon^{-3})$ (Theorem \ref{thm:bigtheorem}). This rate implies the optimal complexity of $O(\epsilon^{-4})$ for finding $\epsilon$-stationary points of smooth objectives \cite{arjevani2023lower,lu2021optimal}. 

\item For nonsmooth stochastic optimization with first-order oracles, ME-DOL achieves the same complexity, $O(\delta^{-1}\epsilon^{-3})$, matching its centralized counterpart (Theorem \ref{thm:ns_fo}). To the best of our knowledge, this is the {\it first} finite-time guarantee for decentralized nonsmooth nonconvex stochastic optimization in the first-order oracle setting. Prior to our work, \citet{Chen_Distributed_weaklyconvex} provided finite-time guarantees only on the Moreau Envelope of {\it weakly-convex} functions.

\item For the zero-order oracle setting, ME-DOL achieves the best known complexity result in terms of $\delta$ and $\epsilon$, i.e.,  $O(\delta^{-1}\epsilon^{-3})$ (Theorem \ref{thm:ns_zo}), which also matches its centralized zero-order counterpart \cite{lin2022gradient}. In the decentralized setting, this rate was previously achieved only with  the variance reduction mechanism \cite{lin2023decentralized}.
\end{itemize}

\subsection{Highlights of Technical Analysis}

\textbf{Randomized Smoothing.} Finite-time analysis of nonsmooth objectives is mainly based on smooth approximations of these objectives. Randomized Smoothing (RS) and Moreau Envelope (ME) are the most common approximation methods. In ME, the original function is approximated with $f_\mu^{ME} (x) = \min_{y\in \mathbb{R}^d} \big\{f(y) + \frac{1}{2\mu} \|y-x\|^2 \big\}$ \citep{davis2019stochastic,scaman2018optimal,ScamanFasterSmoothing}. ME approximation provides theoretical guarantees when applied on structured objectives with regularizer or used with weak-convexity assumption \cite{davis2019stochastic}, but it might be practically unrealistic in some applications that use ReLU neural networks and $\rho$-margin SVMs \cite{tian2022finitetime}.  
On the other hand, in RS the original function $f$ is approximated by $f_{\delta}^{RS}(x) = \mathbb{E}[f(x+\delta u)]$, where $u$ comes from a Gaussian distribution or a uniform distribution on the unit ball. In RS, the smoothness parameter depends on the ambient dimension as $\sqrt{d}$, but it provides favorable theoretical properties (see Proposition \ref{prop:ran_smo}). Specifically, finding a \goldstat point of a nonsmooth $L$-Lipschitz function $f$ can be pursued via finding an $\epsilon$-stationary point of $f_\delta$ with a $\delta L$ approximation error. RS does not require additional assumptions beyond Lipschitz continuity of the function, and this makes RS more tractable in practical applications. In this work, we develop our algorithm based on RS.

\textbf{Nonconvex to Online Conversion.} Another important technique in nonsmooth optimization is a reduction from nonsmooth nonconvex optimization to online learning \cite{cutkosky2023optimal}. In its original form, the method works on centralized optimization, where an online algorithm runs for a certain period, a candidate point is generated, and at the end of the period the online algorithm is restarted. The implementation of the online algorithm can be written explicitly in the form of gradient clipping \cite{kornowski2023algorithm}. From a technical perspective, the use of regret bounds in online learning streamlines the complexity analysis in the nonsmooth optimization. In our paper, we utilize decentralized online algorithm of \citet{ShahDecMD} to address decentralized nonsmooth nonconvex optimization. Compared to the centralized problem \cite{cutkosky2023optimal}, in the decentralized setting, the discrepancy between local variables and global variables makes the analysis more challenging.

\textbf{Geometric Lemma of \cite{kornowski2023algorithm}.} The optimality criterion for nonsmooth analysis is constructed on the Goldstein subdifferential set. The technical result of \citet{kornowski2023algorithm}, which links the Goldstein subdifferential set of $f_{\delta}$ to $f_{\delta+\mu}$ for $\mu,\delta>0$, plays an important role in our analysis. Basically, for the goal of finding a \goldstat point of $f$, we can use a proportion of $\delta$, namely $a\delta$ (for $0 < a < 1$), for smoothing and use the rest of the budget to identify a ($(1-a)\delta,\epsilon$)-stationary point of the smoothed function $f_{a\delta}$ with the smoothness parameter $L_1=O(a^{-1}\delta^{-1})$. This approach allows us to efficiently control the discrepancy terms. 
\begin{remark}
In decentralized nonsmooth nonconvex optimization, our goal is to find a \goldstat point of global function $f$ with randomized smoothing using partial information. For smooth objectives, the complexity of finding $\epsilon$-stationary points  in terms of the smoothness parameter $L_1$ and $\epsilon$ is $O(L_1 \epsilon^{-4})$ \cite{lu2021optimal}. For the nonsmooth objectives, a straightforward application of randomized smoothing with $L_1 = O(\delta^{-1})$ leads to the overall complexity of $O(\delta^{-1} \epsilon^{-4})$, which is sub-optimal. To improve this rate, we develop a technique inspired by \citet{cutkosky2023optimal} and based on decentralized online learning, and we obtain a finite-time bound, in which some of the terms depend on $L_1$. With the geometric lemma of \citet{kornowski2023algorithm} we control the complexity of $L_1$-dependent terms. As a result, we obtain the same complexity rate (up to constant factors) for decentralized nonsmooth nonconvex stochastic optimization as previously achieved in the centralized setting \cite{cutkosky2023optimal}.
\end{remark}

\subsection{Literature Review}
Nonconvex optimization is well-studied under Lipschitz-smoothness assumption. In this setting, the goal is to find an $\epsilon$-stationary point $x$ satisfying $\|\nabla f(x)\| \le \epsilon$. For the deterministic setting, it is well-known that gradient descent (GD) achieves a $O(\epsilon^{-2})$ sample complexity, and this rate is optimal \cite{carmon2020lower}. In the stochastic setting, SGD achieves $O(\epsilon^{-4})$ rate with the assumption of unbiased, bounded variance gradients \cite{ghadimi2013stochastic}. This rate is also optimal as shown by \citet{arjevani2023lower}. We now discuss several strands of related literature.

{\bf Nonsmooth Nonconvex Optimization.} The first non-asymptotic analysis for nonsmooth nonconvex objectives was provided by \citet{zhang2020complexity}, showing that finding $\epsilon$-stationary points in finite time is impossible. Furthermore, it was proved by \citet{kornowski2022oracle} that obtaining a near $\epsilon$-stationary point is also impossible in nonsmooth optimization. Therefore, the goal of finding \goldstat points considered in \citet{zhang2020complexity} is a tractable optimality criterion for nonsmooth objectives. 

{\bf{First-Order Nonsmooth Nonconvex Setting.}} ($\delta,\epsilon$)-Goldstein stationarity of nonsmooth objectives has been analyzed in various settings \citep{tian2022finitetime}. %For stochastic optimization, \citet{davis2019stochastic} showed $O(\epsilon^{-4})$ rate for weakly convex nonsmooth objectives using ME.  
\citet{davis2022gradient} showed that $\tilde{O}(\delta^{-1}\epsilon^{-3} )$ can be achieved for Lipschitz continuous objectives when function values and gradients can be evaluated at points of differentiability. More recently, \citet{cutkosky2023optimal} proved that the $O(\delta^{-1}\epsilon^{-3})$ sample complexity is optimal in the stochastic first-order setting.

{\bf{Zero-Order Nonsmooth Nonconvex Setting.} } Another line of work focuses on zero-order setting in nonsmooth nonconvex optimization. \citet{lin2022gradient} proposed a gradient-free method GFM and its stochastic counterpart SGFM, which achieve $O(d^{\frac{3}{2}}\delta^{-1}\epsilon^{-4})$ sample complexity. \citet{chen2023faster} improved this complexity to  $O(d^{\frac{3}{2}}\delta^{-1}\epsilon^{-3})$ by applying variance reduction. Furthermore, \citet{kornowski2023algorithm} improved the dimension dependence to $O(d \delta^{-1}\epsilon^{-3})$ based on online-to-nonconvex conversion technique introduced by \citet{cutkosky2023optimal}.

{\bf{Deterministic Nonsmooth Nonconvex Setting.}} In the deterministic setting, even in the absence of noise, it is hard to deal with nonsmooth objectives, and randomization is necessary to obtain a dimension independent guarantee. Furthermore, deterministic algorithms require zero-order oracle for finite-time convergence guarantees 
\citep{jordan2022complexity,jordan2023deterministic,tian2022finitetime}.

{\bf{Distributed Smooth Setting.}} In the literature the term distributed may refer to different layers of optimization,  i.e. application, protocol or network topology \cite{lu2021optimal}. In federated learning, it refers to the application layer where each agent uses its local data with shared parameters \cite{mcmahan2016federated}. In fully decentralized scenarios, each agent updates its local parameter using local data and communicates through a connected network. For nonconvex objectives, decentralized SGD has gained a lot of attention \cite{lian2017can} due to the linear speed-up property. Many works analyzed decentralized algorithms under identically distributed data or bounded outer variance assumption for smooth problems \cite{tang2018d,koloskova19a,li2020federated,wang2020tackling,xu2023asynchronous}.

{\bf{Decentralized Nonsmooth Setting.}} For decentralized nonconvex nonsmooth optimization, though the asymptotic analysis was previously explored in \citet{swenson2022distributed}, there exists a scant literature on the {\it finite-time} analysis. For $\lambda$-weakly-convex nonsmooth objectives, \citet{Chen_Distributed_weaklyconvex} provided finite-time guarantees on ME. More recently, \citet{lin2023decentralized} proposed an algorithm (DGFM) that achieves $O(d^{3/2}\delta^{-1}\epsilon^{-4})$ complexity rate in the zero-order setting. In the same setup, they also proposed DGFM+ by incorporating variance reduction to obtain the $O(d^{3/2}\delta^{-1}\epsilon^{-3})$ complexity rate. 

%Due to limited research in decentralized nonsmooth nonconvex stochastic optimization, 
We also focus on decentralized nonsmooth nonconvex stochastic optimization in the present work. We develop a fully decentralized method that mimics a restarting decentralized online learning algorithm. Under mild technical assumptions (e.g., Lipschitz continuity of the objective function and unbiased, bounded variance gradients), we analyze the finite-time performance of the algorithm. We study three settings: (i) smooth first-order, (ii) nonsmooth first-order, and (iii) nonsmooth zero-order. For all of them, we establish the optimal sample complexity as previously derived for the centralized stochastic optimization (Tables \ref{sample-table}-\ref{decentralized-table}). %Most importantly, to the best of our knowledge, we provide the first finite-time characterization of convergence rate in the nonsmooth first-order setting (without weak-convexity assumption \cite{Chen_Distributed_weaklyconvex}), and our result on nonsmooth zero-order setting is without using variance reduction (Table \ref{decentralized-table}).

\begin{table}[t]
\caption{The sample complexity of finding \goldstat points in centralized nonsmooth nonconvex stochastic optimization. $d$ : ambient dimension, $\gamma = f(x_0)-\inf_x f(x)$ where $x_0$ is the initial point, $G^2$: bound on the second moment of stochastic gradient (for first-order methods) and bound on the second moment of Lipschitz constant (for zero-order methods). `*': Oracle requires an additional constraint on the directional derivative. The dependence to dimension $d$ is only reported for zero-order methods.}
\label{sample-table}
\vskip 0.15in
\begin{center}
\begin{small}
\begin{sc}
\begin{tabular}{p{0.8cm} p{1cm} p{2cm} c}
\toprule
Oracle & Method & Reference & Complexity \\
\midrule
First* & SINGD & \cite{zhang2020complexity} 
& $\tilde{O}(\frac{\gamma G^3}{\delta \epsilon^4})$
\\
First & PSINGD & \cite{tian2022finitetime} 
& $\tilde{O}(\frac{\gamma G^3}{\delta \epsilon^4})$
\\
First & O2NC & \cite{cutkosky2023optimal} 
& $O(\frac{\gamma G}{\delta \epsilon^3})$
\\
Zero & SGFM & \cite{lin2022gradient} 
& $O(d^{\frac{3}{2}}(\frac{G^4}{\epsilon^4}+\frac{\gamma G^3}{\delta \epsilon^4}))$
\\

Zero & GFM+ & \cite{chen2023faster} 
& $O(d^{\frac{3}{2}}(\frac{G^3}{\epsilon^3}+\frac{\gamma G^2}{\delta \epsilon^3}))$
\\
Zero & OSNNO  & \cite{kornowski2023algorithm} 
& $O(\frac{d\gamma G^2}{\delta \epsilon^3})$
\\

\bottomrule
\end{tabular}
\end{sc}
\end{small}
\end{center}
\vskip -0.1in
\end{table}

\begin{table}[t]
\caption{The sample complexity of finding \goldstat points for decentralized nonsmooth nonconvex stochastic optimization. `*' : weakly convex setting, `**' : variance reduction. The dependence to dimension $d$ is only reported for zero-order methods.}
\label{decentralized-table}
\vskip 0.15in
\begin{center}
\begin{small}
\begin{sc}
\begin{tabular}{p{0.8cm} p{1.45cm} p{2.749cm} p{1.5cm}}
\toprule
Oracle & Method & Reference & Complexity \\
\midrule
First* & DPSM & \cite{Chen_Distributed_weaklyconvex} 
& $O(\epsilon^{-4})$
\\
Zero & DGFM & \cite{lin2023decentralized} 
& $O(d^{\frac{3}{2}}\delta^{-1}\epsilon^{-4})$
\\
Zero** & DGFM+ & \cite{lin2023decentralized} 
& $O(d^{\frac{3}{2}}  \delta^{-1} \epsilon^{-3})$
\\
Zero & ME-DOL & {\bf Our Work}
& $O(d\delta^{-1} \epsilon^{-3})$
\\
First & ME-DOL & {\bf Our Work}
& $O(\delta^{-1} \epsilon^{-3})$
\\

\bottomrule
\end{tabular}
\end{sc}
\end{small}
\end{center}
\vskip -0.1in
\end{table}

\section{Problem Setting}
{\bf Notation:} We denote by $\norm{x}$ the Euclidean norm, by $[n]$ the set $\{1,2,3,...,n\}$, by $B(x,\delta):= \{y \in \mathbb{R}^d : \norm{y-x}\le\delta\}$, by $conv(\cdot)$ the convex hull operator, and by $\un(A)$ the uniform measure over a set $A$. $\norm{\cdot}_F$ denotes the Frobenius norm and $\norm{\cdot}_2$ denotes the spectral norm. We use the standard notation $O(\cdot)$, $\Theta(\cdot)$ , $\Omega(\cdot)$ to hide the absolute constants and $\tilde{O}(\cdot)$ to hide poly-logarithmic factors. $1_d$ and $\mathbb{S}^{d-1}$ denote the vector of all ones and the unit sphere in $\mathbb{R}^d$, respectively.

{\bf Network Setup:} In decentralized learning, we have $n$ agents that communicate through a network. We assume that the network is connected, i.e., there exists a (potentially multi-hop) path from any agent $i\in [n]$ to $j\neq i$. The network topology is governed by a symmetric doubly stochastic matrix $P=[P_{ij}]_{i,j=1}^n \in \mathbb{R}^{n\times n}$, where $1_n^\top P = 1_n^\top$ and $P 1_n = 1_n$ (Assumption \ref{assumption:dsm}). Throughout the learning process, agent $i$ receives only information about its local function $f^i$ in the form of stochastic gradients or noisy function evaluations. Note that $P_{ij}\in [0,1]$ and if $P_{ij}=0$ agents $i$ and $j$ do not directly share information with each other. However, if $P_{ij}>0$ agents share their decision variables as described in the ME-DOL (Algorithm \ref{alg:decen}). As such, the neighborhood of agent $i$ is defined as $\mathcal{N}_i:=\{j\in [n]: P_{ij}>0\}$.

\textbf{Information Oracles:} We assume that each agent $i\in [n]$ has access to its information oracles. In the first-order setting, the oracle $\oracle{f}$ returns stochastic gradient at query point $x$ given by $\oracle{f}^{i}(x) = \nabla F^i (x,\xi_i)$. For the first-order oracles we assume that the oracle returns unbiased estimates of the gradient with bounded variance (Assumption \ref{assumption:gradientoracle}). In the zero-order setting, agents have access to the stochastic function value oracle $\oracle{z}$ at query point $x$, that is $\oracle{z}^{i}(x) =  F^i (x,\xi_i)$ with Assumption \ref{assumption:Lipschitz} in place.

\subsection{Stationarity Metric in Nonsmooth Analysis}
In nonconvex smooth optimization problems, finding an $\epsilon$-stationary point $x$, i.e. $\|\nabla f(x) \| \le \epsilon$ is a well-known tractable optimality condition.  For nonsmooth objectives, a more relaxed criterion called near $\epsilon$-stationarity can be considered for a point $x$ with $\min\{ \norm{g} : g \in \cup_{y \in B(x,\delta)} \partial f(y) \}\le \epsilon$. However, both could be intractable criteria for nonsmooth objectives \cite{kornowski2022oracle}. By Rademacher's Theorem, Lipschitz continuous functions are almost everywhere differentiable. For this class of functions, we can study ($\delta,\epsilon$)-stationarity. Let us first define Goldstein $\delta$-subdifferential as follows.
\begin{definition}
\label{def:goldsubdiff}
Goldstein $\delta$-subdifferential of $f$ at $x$ 
is the set $$\partial_{\delta}f(x) := conv(\cup_{y \in B(x,\delta)} \partial f(y)),$$    where the Clarke subdifferential set $ \displaystyle \partial f(x) := conv\{g : g = \lim_{x_s \rightarrow x} \nabla f(x_s)\}$. 
\end{definition}
Since Goldstein $\delta$-subdifferential is the convex hull of a set of Clarke subdifferentials, it is possible that an element of Goldstein $\delta$-subdifferential is not an element of Clarke subdifferential set of points $y \in B(x,\delta)$ for a nondifferentiable function $f$. An example of a function that has a \goldstat point that is not near $\epsilon$-stationary is given in \citet{kornowski2022oracle} (Proposition 2).
\begin{definition}\label{def:goldstein}
Given a Lipschitz function $ f:\mathbb{R}^d \rightarrow \mathbb{R} $, a point $x \in \mathbb{R}^d$ and $\delta > 0$ , denote  $\gradnormdelta{f}{x} := \min\{\|g\| : g \in \partial_{\delta} f(x) \} $. A point $x$ is called a \goldstat point of $f(\cdot)$ if $\gradnormdelta{f}{x} \le \epsilon$. 
\end{definition}
This is a weaker notion than $\epsilon$-stationarity or near $\epsilon$-stationarity. In case of differentiable functions with $L_1$ Lipschitz gradients, an $(\frac{\epsilon}{3L_1},\frac{\epsilon}{3} )$-Goldstein stationary point is also $\epsilon$-stationary \cite{zhang2020complexity}.

\subsection{Properties of Randomized Smoothing}
In the nonsmooth analysis, finding \goldstat points of the global function $f(x) = \frac{1}{n} \sum_{i=1}^{n}{f^i}(x)$ is a reasonable and tractable optimality criterion using {\it randomized smoothing}.

\begin{definition} 
Given an $L$-Lipschitz function $f$, we denote its smoothed surrogate as $f_{\delta}(x) := \mathbb{E}_{u \thicksim \mathcal{P} }[f(x+\delta u )]$, where $\mathcal{P}$ is the uniform distribution on the unit ball, i.e., \un($B(0,1)$).
\label{def:uniformsmoothing}
\end{definition}

\begin{proposition}\label{prop:ran_smo}\cite{lin2022gradient}
Suppose that the function $f : \mathbb{R}^d \rightarrow \mathbb{R}$ is $L$-Lipschitz. Then, it holds that:
\begin{itemize}
\item $|f_\delta(\cdot) - f(\cdot)| \le \delta L$.
\item $f_\delta(\cdot)$ is $L$-Lipschitz.
\item $f_\delta(\cdot)$ is differentiable with $c\sqrt{d}L \delta^{-1}$-Lipschitz gradients for a numeric constant $c>0$.
\item $\nabla f_\delta(\cdot) \in \partial_\delta f(\cdot)$ , where $\partial_\delta f(\cdot)$ is the Goldstein subdifferential.
\end{itemize}
\end{proposition}
RS allows us to work with a ``smoothed" objective $f_{\delta}$ with the cost of $\delta L$ approximation error. Also, the last item in Proposition \ref{prop:ran_smo} links finding a \goldstat point of a nonsmooth objective $f$ with finding an $\epsilon$-stationary point of the smoothed objective $f_{\delta}$. Furthermore, we will use the following lemma for nonsmooth analysis.
\begin{lemma} \cite{kornowski2023algorithm}
For any $\delta,\mu \ge 0 : \partial_{\mu} f_{\delta}(x) \subseteq \partial_{\mu+\delta}f(x)$ . 
\label{lem:geo_gold}
\end{lemma}
\begin{proposition}
\label{prop:deltanorm}
By Lemma \ref{lem:geo_gold} and Definition \ref{def:goldstein} for $\norm{\cdot}_{\delta}$ we have $\norm{\nabla f(x)}_{\delta} \le \norm{\nabla f_{a\delta}(x)}_{(1-a)\delta} $ for any $a \in (0,1)$. 
\end{proposition}
For the task of finding \goldstat points of a function, $\delta$ denotes the radius of the ball as in Definition \ref{def:goldsubdiff}. Lemma \ref{lem:geo_gold}
allows us to distribute $\delta$ such that we can use a portion of it for smoothing and allocate the remaining part to the radius of the ball around the critical point. In Proposition \ref{prop:deltanorm} by choosing $a=\frac{1}{2}$ we have $\norm{\nabla f(x)}_{\delta} \le \|\nabla f_{\frac{\delta}{2}}(x)\|_{\frac{\delta}{2}}$.

Using these lemmas and randomized smoothing, we can facilitate our convergence analysis. In the context of decentralized optimization, we can use the surrogate function $f_\delta$ of global function $f$, and by the linearity of expectation we have $f_{\delta} = \frac{1}{n} \sum_{i=1}^{n}{(f^i)_{\delta}}$. Furthermore, using the following lemma, summation of gradients of smoothed functions  $\frac{1}{n} \sum_{i=1}^{n}\nabla (f^i)_{\delta}(x_i)$ can be related to $\frac{1}{n} \sum_{i=1}^{n}\nabla f_{\delta}(x_i)$. %and then by applying Lemma \ref{lem:geo_gold} we can extract information about $\partial_{\mu + \delta} f(\bar{x})$.

\begin{lemma}
\label{lem:decentralized}
Suppose that  $n$ local functions $\{f^i\}_{i=1}^{n}$ have $L_1$-Lipschitz gradients and $f(x) = \frac{1}{n} \sum_{i=1}^n f^i(x)$. Consider the set of points $\{w_{t,i}\}$ for $i\in [n], t\in [T]$, and let $\bar{w}_t = \frac{1}{n} \sum_{i=1}^{n} w_{t,i} $ and $\norm{w_{t,i} -\bar{w}_t} \le r, \forall i \in [n], \forall t \in [T] $. Then, we have {\small
$$\norm{\frac{1}{nT} \finitesum{t}{T} \finitesum{i}{n} \nabla f(w_{t,i}) } \le \norm{\frac{1}{nT} \finitesum{t}{T} \finitesum{i}{n} \nabla f^i(w_{t,i}) } +2rL_1.$$ }
\end{lemma}

\subsection{Assumptions}
We assume that agents communicate synchronously through the network. For example, agent $i$ takes a weighted average of the decision variables in its neighborhood as follows
$$ \displaystyle y_{t,i} = \sum_{j \in \mathcal{N}_i }P_{ij}x_{t,j}=\sum_{j=1}^nP_{ij}x_{t,j} , \forall i \in [n],$$
as elaborated in Algorithm \ref{alg:decen}. The communication matrix $P$ is fixed over time and satisfies the following assumption. 
\begin{assumption}
The network is connected and the communication matrix $P \in \mathbb{R}^{n \times n}$ is symmetric and doubly stochastic. $\rho$ denotes the second largest singular value of the matrix $P$. Given that the network is connected, we have that $\rho \in [0,1)$. 
\label{assumption:dsm}
\end{assumption}
Assumption \ref{assumption:dsm} is widely used in the decentralized optimization literature (see e.g., \cite{ShahDecMD}). Regardless of whether the network structure is fixed or time-varying, some connectivity assumption is needed to solve the global problem. Here, the quantity $\rho$ determines the connectivity of the network, and a smaller $\rho$ indicates a more well-connected network topology.

\begin{assumption}
We assume that local objective functions have the form $f^i (x) = \mathbb{E}_{\xi} [F^i (x,\xi)] $, where $\xi$ denotes the random index. The stochastic component of local functions $F^i(\cdot,\xi) : \mathbb{R}^d \rightarrow \mathbb{R}$ is $L(\xi)$-Lipschitz for any $\xi$, i.e., it holds that $$\left|F^i(x,\xi)-F^i(y,\xi)\right| \le L(\xi) \norm{x-y},$$ for any $x,y \in \mathbb{R}^d$ and $i\in [n]$. $L(\xi)$ has a bounded second moment such that $\mathbb{E}_{\xi} [L(\xi)^2] \le L^2$. %We further assume that $f^i$ is directionally differentiable.
\label{assumption:Lipschitz}
\end{assumption}
We note that Assumption \ref{assumption:Lipschitz} is weaker than assuming that $F^i(\cdot,\xi)$ is $L$-Lipschitz \cite{chen2023faster,kornowski2023algorithm}. Taking expectation from above, it can be shown that local functions $f^i$ are Lipschitz continuous. However, we do {\it not} assume that gradients are Lipschitz. Furthermore, directional differentiability holds for commonly used nonsmooth functions (e.g., ReLU) and enables the use of Lebesgue path integrals \cite{zhang2020complexity}.

\begin{assumption}
The local objectives $f^i:\mathbb{R}^d \rightarrow \mathbb{R}$ are lower bounded $(f^i)^* := \inf_x f^i(x) > -\infty$. Therefore, the global function $f$ is also lower bounded, and we define $\gamma$ such that $f(\bar{x}_0) - \inf_x f(x) \le \gamma $, where $\bar{x}_0$ is the average of initial points (among agents) for the algorithm.
\label{assumption:suboptimality}
\end{assumption}

We also make the following standard assumptions on the stochastic gradients \cite{ShahDecMD,zhang2020complexity} .
\begin{assumption}\label{assumption:gradientoracle}
We assume that the first-order oracle returns unbiased, bounded variance estimate of the gradient such that $\mathbb{E}[\nabla F^i(x,\xi)] = \nabla f^i(x)$ and $\probexp{ \norm{ \nabla F^i(x,\xi) - \nabla f^i(x) }^2}  \le \sigma^2 $. Furthermore, we assume that the second moment of the stochastic gradient is bounded such that  $\probexp{\norm{\nabla F^i(x,\xi)}^2} \le G^2 $.
\end{assumption}

\section{Algorithm and Main Technical Results}

In this section, we present our decentralized algorithm for finding a \goldstat point of the global objective $f$ in \eqref{eq:global}. Our algorithm is termed Multi Epoch Decentralized Online Learning (ME-DOL), for which we establish the sample complexity in different settings. First, we present our result on smooth nonconvex objectives (Theorem \ref{thm:bigtheorem}), and we then extend our results to nonsmooth nonconvex objectives with randomized smoothing in first-order and zero-order settings (Theorems \ref{thm:ns_fo} and \ref{thm:ns_zo}, respectively). Our technique is based on a reduction from nonsmooth nonconvex decentralized optimization to decentralized online learning.

\begin{algorithm}[tb]
   \caption{Multi Epoch Decentralized Online Learning}
   \label{alg:decen}
\begin{algorithmic}
   \STATE {\bfseries Input:} $\delta' \in \mathbb{R}_{\geq0}$, $K \in \mathbb{N}$ , $T \in \mathbb{N}$, decentralized online learning algorithm $\mathcal{A}$ with bounded domain $\mathcal{D}$, doubly stochastic communication matrix $P$.
   \STATE {\bf Initialize:} $y_{T,i}^0 = 0$ for all $i \in [n]$.
   \FOR{$k=1$ {\bfseries to} $K$ } 
   \STATE Restart $\mathcal{A}$ 
   \STATE Let $y_{0,i}^k = y_{T,i}^{k-1}$ ,  $\forall i \in [n]$ 
   \FOR{$t=1$ {\bfseries to} $T$,  $\forall i \in [n]$ }
   \STATE Get $\Delta_{t,i}^k$ for all agents from $\mathcal{A}$ (Algorithm \ref{alg:DOL})
   \STATE $x_{t,i}^k = y_{t-1,i}^k+\Delta_{t,i}^k$ 
   \STATE $s_{t,i}^k\sim \un[0,1]$
   \STATE $w_{t,i}^k = y_{t-1,i}^k+s_{t,i}^k\Delta_{t,i}^k$ 
   \STATE $y_{t,i}^k = \sum_{j=1}^n P_{ij} x_{t,j}^k= \sum_{j\in\mathcal{N}_i} P_{ij} x_{t,j}^k  $

   \IF{Information Oracle == Zero-Order}
   \STATE $g_{t,i}^k$ = Zero-Order Gradient($F^i,w_{t,i}^k,\delta ' , \xi_{t,i}^k)$
   \ELSE
   \STATE $g_{t,i}^{k}$ = First-Order Gradient($F^i,w_{t,i}^k,\delta ' , \xi_{t,i}^k)$

   \ENDIF
   \STATE Send $g_{t,i}^k$ to $\mathcal{A}$ as gradient
   \ENDFOR
   \STATE Set $\bar{w}^k = \frac{1}{nT}\sum_{i=1}^{n} \sum_{t=1}^{T} w_{t,i}^k $ for $k\in [K]$
   \ENDFOR
   \STATE {Sample $w^{out} \sim \un\{ \bar{w}^1,..., \bar{w}^K\} $}
   \STATE{\bfseries Output : } $w^{out}$
\end{algorithmic}
\end{algorithm}

In Algorithm \ref{alg:decen}, we have periods of length $T$, where in each epoch $k\in [K]$ a decentralized online algorithm $\mathcal{A}$ is used to generate action $\Delta^k_{t,i}$ for agent $i\in [n]$ at iteration $t\in [T]$. The action space $\mathcal{D}$ is bounded such that $\norm{\Delta} \le D$ for all $\Delta \in \mathcal{D}$. At each epoch $k$, $\mathcal{A}$ must run on a linear optimization problem with the objective of agent $i$ as $\sum_{t=1}^T\innerproduct{\Delta,g^k_{t,i}}$, where $g_{t,i}^k$ is the stochastic gradient (respectively, approximation of the stochastic gradient using noisy function evaluations) in the first-order (respectively, zero-order) setting. We use Algorithm \ref{alg:DOL} \cite{ShahDecMD} for $\mathcal{A}$. Based on action $\Delta^k_{t,i}$ the variable $x_{t,i}^k$ is generated and then averaged over neighborhood of $i$ to get $y_{t,i}^k$. The set of $nT$ points $w^k_{t,i}$ $\forall i\in[n], \forall t\in [T]$, at which the gradients are evaluated, are averaged to produce the final output of each epoch, denoted as $\bar{w}^k$. These points are proposed as candidates for identifying a \goldstat point of the global function $f$. Algorithm \ref{alg:decen} outputs a randomly selected candidate point $\bar{w}^{l}$ where $l\sim \un[K]$.
\begin{remark}
Note that $P_{ij}=0$ if agents $i$ and $j$ are not neighbors. Therefore, updates $y_{t,i}^k = \sum_{j=1}^n P_{ij} x_{t,j}^k$ (in Algorithm \ref{alg:decen}) and $\Delta^k_{t+\frac{1}{2},i} = \finitesum{j}{n} P_{ij}\Delta^k_{t,j}$ (in Algorithm \ref{alg:DOL}) do not contradict the decentralized nature of the learning. Basically, for each agent $i\in [n]$, the sum always reduces to a weighted averaging over the neighborhood of agent $i$.
\end{remark}

\begin{remark}
\label{rem:DOLMD }
The original implementation of Algorithm \ref{alg:DOL} in \citet{ShahDecMD} is based on mirror descent, but here we use the Euclidean distance as the generator of Bregman divergence, reducing the algorithm to decentralized online gradient descent.  
\end{remark}

\begin{algorithm}[tb]
    \caption{First-Order Gradient($F, x, \delta', \xi$)}
   \label{alg:firstorder}
\begin{algorithmic}
    \STATE {\bfseries Input:} Function $F$, point $x$, smoothing parameter $\delta '$, random seed $\xi$.
    \STATE Sample $z \thicksim $  \un($B(0,1)$) 
    \STATE $ g = \nabla F( x + \delta ' z, \xi ) $
    \STATE{\bfseries Output:} $g$

\end{algorithmic}
\end{algorithm}

\begin{algorithm}[tb]
    \caption{Zero-Order Gradient($F, x,\delta' , \xi$)}
   \label{alg:zeroorder}
\begin{algorithmic}
    \STATE {\bfseries Input:} Function $F$, point $x$, smoothing parameter $\delta '$, random seed $\xi$, dimension $d$.

    \STATE Sample $ z \thicksim $ \un($\mathbb{S}^{d-1}$) 
    
    \STATE Evaluate $F(x + \delta ' z, \xi )$ and $F(x - \delta ' z, \xi )$ 

    \STATE $g = \frac{d}{2\delta '} \Big(F(x + \delta 'z, \xi )-F(x - \delta 'z, \xi )\Big) z $

    \STATE{\bfseries Output:} $g$

\end{algorithmic}
\end{algorithm}

\begin{algorithm}[tb]
    \caption{ Decentralized Online Optimization Algorithm $\mathcal{A}$ \cite{ShahDecMD} }
   \label{alg:DOL}
\begin{algorithmic}
    \STATE {\bfseries Input:} Domain $\mathcal{D}$, doubly stochastic matrix $P$, learning rate $\eta$, stochastic gradients $g_{t,i}^k$.  \STATE {\bfseries Initialize:} $\Delta_{\frac{1}{2},i}^k = 0$ for all $i\in [n]$ and $k\in[K]$.
    \STATE {\bfseries Iterations : } Step $t \ge 1$, update for each $i \in [n]$:
        \STATE $\Delta_{t,i}^k = \underset{\Delta \in \mathcal{D}}{\text{argmin}} \Big\{\eta \innerproduct{\Delta,g_{t-1,i}^k} + \frac{1}{2} \big\|\Delta - \Delta_{t-\frac{1}{2},i}^k \big\|^2 \Big\} $
        \vspace{.15cm}
        \STATE $\Delta_{t+\frac{1}{2},i}^k = \finitesum{j}{n} P_{ij}\Delta_{t,j}^k $   
\end{algorithmic}
\end{algorithm}

\subsection{Challenges in the Analysis of Decentralized Algorithm}
In centralized optimization, the difference of function values in consecutive iterations, $f(x_t) - f(x_{t-1})$, depends on the update rule. The update rule can be written as $x_t = x_{t-1}+\Delta_t$. For example, in SGD,  $\Delta_t = -\eta \nabla F (x_{t-1},\xi)$. Various algorithms use the past information to generate $\Delta_t$ or process the latest information as in the normalized gradient descent \cite{murray2019revisiting} or gradient clipping \citep{zhang2019gradient}. For any algorithm with the update rule $x_t = x_{t-1}+\Delta_t$, one can write $f(x_t) = f(x_{t-1}) + \langle \Delta_t , \nabla_t \rangle$ where $\nabla_t = \int_{0}^{1}\nabla f(x_{t-1}+s\Delta_{t})ds$. \citet{cutkosky2023optimal} observed that $\Delta_t$ can be generated via an online learning algorithm (e.g., online gradient descent (OGD)), in which case the summation of the differences for $T$ rounds can be written as follows 
\begin{align*}
     &f(x_T) - f(x_0) = \sum_{t=1}^{T}  \innerproduct{\Delta_t , \nabla_t}\\
&~~~~~= \sum_{t=1}^{T}  \innerproduct{ g_t , \Delta_t - u } + \sum_{t=1}^{T}  \innerproduct{\nabla_t - g_t , \Delta_t } + \sum_{t=1}^{T}  \innerproduct{ g_t , u }, 
\end{align*}   
where $\{g_t\}$ are the stochastic gradients provided to the online learning algorithm. This equation holds for any $u$ in hindsight, and the first summation term corresponds to the regret of the online algorithm. The second term has expectation equal to $0$ (if gradients are unbiased), and for the last term we have freedom to select an optimal $u$. 

In the decentralized counterpart of this conversion, which is the focal point of our analysis, the main difference is that the change in the global function depends on the average action, i.e., $\bar{\Delta}_t = \bar{x}_t - \bar{x}_{t-1}$, but $f(\bar{x}_t) = f(\bar{x}_{t-1}) + \langle \bar{\Delta}_t , \tilde{\nabla}_t \rangle$ where $\tilde{\nabla}_t = \int_{0}^{1}\nabla f(\bar{x}_{t-1}+s\bar{\Delta}_{t})ds $. The key technical challenge is the discrepancy between $\tilde{\nabla}_t $ and $\bar{\nabla}_t = \sum_{i=1}^n{\nabla_{t,i}/n}$, where $\nabla_{t,i} := \int_{0}^{1}\nabla f^{i}(y_{t-1,i}+s\Delta_{t,i})ds$. In the decentralized analysis, we have the decentralized regret term and  an additional discrepancy term that arises due to the difference between $\tilde{\nabla}_t- \bar{\nabla}_t$, which requires careful analysis.

\subsection{Smooth Analysis}
Let us now present our first result on smooth objectives in the following theorem.
\begin{theorem}
\label{thm:bigtheorem}
Let Assumptions \ref{assumption:dsm}, \ref{assumption:suboptimality}, \ref{assumption:gradientoracle} hold and further assume that $f^i$ is $L$-Lipschitz and has $L_1$-Lipschitz gradient for all $i\in [n]$. Let $\delta,\epsilon \in (0,1)$ and choose 
\begin{itemize}
\item $N := KT = \Theta(\delta^{-1} \epsilon^{-3}(1-\rho)^{-2})$,
\item {$T = \Theta((1-\rho)^{\frac{1}{3}}\left(\delta N \right)^{\frac{2}{3}})$,}
\item $D=\frac{\delta(1-\rho)}{2T\sqrt{n}}$,
\item {$\eta =\Theta(\sqrt{(1-\rho)}\frac{D}{\sqrt{T}})$.}
\end{itemize}
Then, running Algorithm \ref{alg:decen} with $\delta'=0$ for $N$ rounds gives an output that satisfies the following inequality for the global function $f(x)=\frac{1}{n}\sum f^i(x)$,
$$\mathbb{E}_{k\sim \un[K]}\left[ \gradnormdelta{f}{\bar{w}^k}\right] %O\Big((\delta N)^{-\frac{1}{3}}(1-\rho)^{-\frac{2}{3}}\Big)
\leq \epsilon.$$
\end{theorem}

\begin{remark}
\label{rem:one}
For smooth objectives this result implies that a \goldstat point can be found in  $ N =O(\delta^{-1}\epsilon^{-3})$ iterations. This rate results in the optimal complexity of $O(\epsilon^{-4})$ for finding an $\epsilon$-stationary point of nonconvex smooth objectives \cite{arjevani2023lower,lu2021optimal} as $\delta=O(\epsilon)$. Furthermore, the dependence of $N$ to $(1-\rho)^{-2}$ indicates that in a well-connected network (smaller $\rho$), we need less iterations to find a \goldstat point.
\end{remark}
For smooth objectives we do not need randomized smoothing, so we choose $\delta'=0$. We can use the full budget $\delta$ for the search radius, so we set $D = \frac{\delta (1-\rho)}{2T\sqrt{n}}$ in order to satisfy $\|\bar{w}^k -w_{t,i}^k \| \le \delta$. The complete proof of Theorem 1 can be found in the Appendix (Section \ref{proof:Thm31}).

%In the proof of Theorem 1 (Equation \ref{eq:smooth_ineq} in the Appendix), the bound has two components with respect to the smoothness parameter $L_1$.  As $N$ increases, the $L_1$-dependent term is dominated by the $L_1$-independent term, which allows us to use RS for nonsmooth objectives.

\subsection{Challenges in Nonsmooth Analysis}
For nonsmooth objectives, we utilize randomized smoothing. The straightforward application of randomized smoothing, such as merely replacing the task of finding a \goldstat point of a nonsmooth function $f$ with the task of finding an $\epsilon$-stationary point of the smoothed function would result in the sub-optimal complexity of $O(\delta^{-1}\epsilon^{-4})$ since the optimal rate for decentralized smooth nonconvex objectives is $O(L_1\epsilon^{-4})$ \cite{lu2021optimal}, and $L_1 = O(\delta^{-1})$ for the smoothed function according to Proposition \ref{prop:ran_smo}.

To address this, we must control $\delta$ that affects the smoothness parameter $L_1$. In the proof of Theorem \ref{thm:bigtheorem}, we have the following inequality (see Equation \ref{eq:smooth_ineq}), which also plays an important role in the nonsmooth analysis.
\begin{align*}
&\epstermglobal\\
&~~~~~~~~\le \frac{2\gamma T \sqrt{n}  }{\delta N (1-\rho)} + \frac{\sigma}{\sqrt{nT}} + \frac{c_1}{\sqrt{T}} + \frac{\delta L_1 (1-\rho) c_3}{2T \sqrt{n}}.
\end{align*}
 However, for nonsmooth objectives, $f$ must be replaced by a smoothed function $f_{(1-a)\delta}$ in the left-hand side, and using Proposition \ref{prop:deltanorm}, we can consider finding an ($a\delta,\epsilon$)-stationary point of $f_{(1-a)\delta}$ for $0<a<1$. A larger $a$ increases the ``radius of possible stationarity" but decreases the ``smoothness", and we need to balance this trade-off. 

\subsection{Nonsmooth Analysis with First-Order Oracle}
For nonsmooth objectives, we can choose $a=\frac{1}{2}$, and the goal is to find a ($\frac{\delta}{2},\epsilon$)-stationary point of $f_{\frac{\delta}{2}}$. Then, we can extend the result of Theorem \ref{thm:bigtheorem} with the smoothness parameter $L_1 = 2c\sqrt{d}L \delta^{-1}  $, where $c$ is a constant that depends on the geometry of the problem (see Appendix \ref{app:constant}).
\begin{theorem}
\label{thm:ns_fo}
Let $\delta,\epsilon \in (0,1)$. Suppose that Assumptions  \ref{assumption:dsm}, \ref{assumption:suboptimality}, \ref{assumption:gradientoracle} hold and that $f^i$ is $L$-Lipschitz for all $i\in [n]$. Choose $N, T, \eta$ as in Theorem \ref{thm:bigtheorem}, and set $D=\frac{\delta(1-\rho)}{4T\sqrt{n}}$. Then, running Algorithm \ref{alg:decen} with $\delta' = \frac{\delta}{2}$ for $N$ rounds gives an output that satisfies the following inequality
$$\mathbb{E}_{k\sim \un[K]}\left[ \gradnormdelta{f}{\bar{w}^k} \right] \le c_8 (\delta N)^{-\frac{1}{3}}\leq \epsilon,$$ 
where $c_8=O((1-\rho)^{-\frac{2}{3}})$ and it does not depend on $\delta$ and $N$ (see Appendix \ref{proof:ns_fo} for the exact quantification of $c_8$).
\end{theorem}
To the best of our knowledge, the above theorem establishes the first complexity rate for nonsmooth nonconvex functions in decentralized stochastic optimization (without weak-convexity assumption). In terms of $\delta$ and $\epsilon$, the rate $N=O(\delta^{-1}\epsilon^{-3})$ matches the best known rates in the centralized setting \cite{cutkosky2023optimal}. This rate also recovers the optimal results in the smooth nonconvex setting when $\delta=O(\epsilon)$. 

\subsection{Nonsmooth Analysis with Zero-Order Oracle}
In the first-order setting, we used Assumption \ref{assumption:gradientoracle} that implies access to unbiased, bounded variance gradient estimates. In the zero-order setting, agents do not have access to the gradient information, but they can estimate the gradient using noisy function values. The output of Algorithm \ref{alg:zeroorder} provides an unbiased gradient estimator with bounded variance \cite{shamir2017optimal}, and it can be used in lieu of stochastic gradients returned by the first-order oracle.
\begin{lemma}\cite{kornowski2023algorithm}
Let $w = x + s\Delta$ be a point with $s\sim \un[0,1]$. The gradient estimator
$$ g = \frac{d}{2\delta '} \Big(F(x + s\Delta + \delta 'z, \xi )-F(x + s\Delta - \delta' z, \xi )\Big)z,$$ as generated by Algorithm \ref{alg:zeroorder} satisfies the following conditions
$$\mathbb{E}_{\xi,z} [g |x, s, \Delta  ] = \nabla f_{\delta '}(x+s\Delta) = \nabla f_{\delta '}(w),$$
and
$$\mathbb{E}_{\xi,z}[\norm{g}^2 |x,s,\Delta ] \le \zeroordersecbound.$$
\label{lemma:zeroorder}
\end{lemma}
The bound on the second moment helps us  replace $G$ and $\sigma$ in the first-order setting (Assumption \ref{assumption:gradientoracle}) by a quantifiable constant. Running ME-DOL using (the noisy version of) smoothed functions $f_{\delta}^i$ and the gradient estimator in Algorithm \ref{alg:zeroorder}, we have the following convergence guarantee for the zero-order setting.
\begin{theorem}
\label{thm:ns_zo}
Let $\delta,\epsilon \in (0,1)$. Suppose that Assumptions \ref{assumption:dsm}, \ref{assumption:Lipschitz}, \ref{assumption:suboptimality} hold and that the zero-order oracle returns unbiased estimates of the function values. Choose $N, T, \eta$ as in Theorem \ref{thm:bigtheorem}, and set $D=\frac{\delta(1-\rho)}{4T\sqrt{n}}$. Then, running Algorithm \ref{alg:decen} with $\delta' = \frac{\delta}{2}$ for $N$ rounds gives an output that satisfies the following inequality
$$\mathbb{E}_{k\sim \un[K]}\left[ \gradnormdelta{f}{\bar{w}^k}  \right] \le c_{11} (\delta N)^{-\frac{1}{3}}\leq \epsilon,$$ 
where $c_{11}=O(d^{\frac{1}{3}}(1-\rho)^{-\frac{2}{3}})$ and it does not depend on $\delta$ and $N$ (see Appendix \ref{proof:ns_zo} for the exact quantification of $c_{11}$).
\end{theorem}
Similar to previous results, the complexity of finding a \goldstat point is $N=O(\delta^{-1}\epsilon^{-3})$ in terms of $\delta$ and $\epsilon$.  In decentralized zero-order nonsmooth nonconvex stochastic optimization, our result matches the best known rate $O(\delta^{-1}\epsilon^{-3})$ as shown in Table \ref{decentralized-table} {\it without} recourse to variance reduction.

\begin{remark}
\label{rem:dimdep}In the zero-order setting (Table \ref{decentralized-table}), the dimension dependence of our algorithm is $O(d)$, which improves upon DGFM and DGFM+, where the dimension dependence is $O(d^{\frac{3}{2}})$. This result also matches with the optimal dimension dependence in the centralized setting following the analysis of \citet{kornowski2023algorithm}.
\end{remark}

\section{Numerical Experiments}
To validate the performance of our algorithm, we conduct experiments on several datasets\footnote{Codes for numerical experiments are available at https://github.com/emreesahinoglu/Decentralized-Nonsmooth.git}. %The feature vectors are normalized before optimization.

{\bf{Model.}}
We consider the nonconvex penalized SVM with capped-$\ell_1$ regularizer. The model trains a binary classifier $x \in \mathbb{R}^d$ on the training data $\{ a_i,b_i\}_{i=1}^m$, where $a_i \in \mathbb{R}^d$ and ${b_i} \in \{-1,1\}$ are the (normalized) feature vector and label for the $i$-th sample, respectively. Local objective functions can be written as
$$f^i(x) = \frac{1}{m_i} \sum_{j=1}^{m_i} l(b_i^j (a_i^j)^\top x) + \nu(x),$$
where $l(y) = \max\{1-y,0\}$, $m= \sum_{i=1}^n m_i$, $\nu(x)=\lambda \sum_{j=1}^d \min\{|x(j)|,\alpha\}$, and $\lambda,\alpha >0$. Similar to experiments of \citet{lin2023decentralized}, we set $\lambda = 10^{-5}/n$ and $\alpha =2$. For each dataset, we divide the training samples equally among agents, i.e., $m_i=m/n$. 

{\bf{Setup.}} We consider a network of $n=20$ agents with a ring topology. Hyper-parameters of our algorithm, $\eta$ and $D$, are selected based on the theorems, where $\eta=\Theta(D/\sqrt{T})$. We set $\eta = 0.01 \times D $ and vary $D$ in the range of $10^{-4}$ to $10^{-2}$ in different experiments. %We also set $\delta=0.001$.

\begin{figure}
%\vskip -0.2in
\begin{center}
\centerline{\includegraphics[width=\columnwidth]{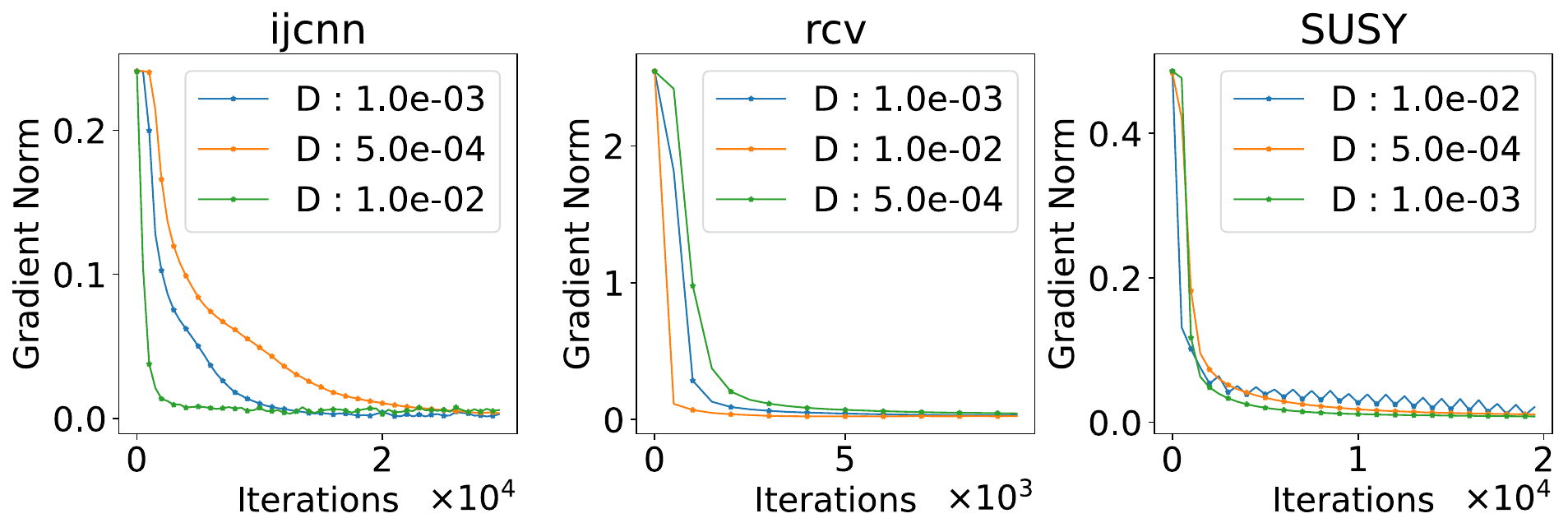}}
\caption{Evaluation of the gradient norm in the first-order setting.}
\label{fig_grad_fo}
\end{center}
\begin{center}
\centerline{\includegraphics[width=\columnwidth]{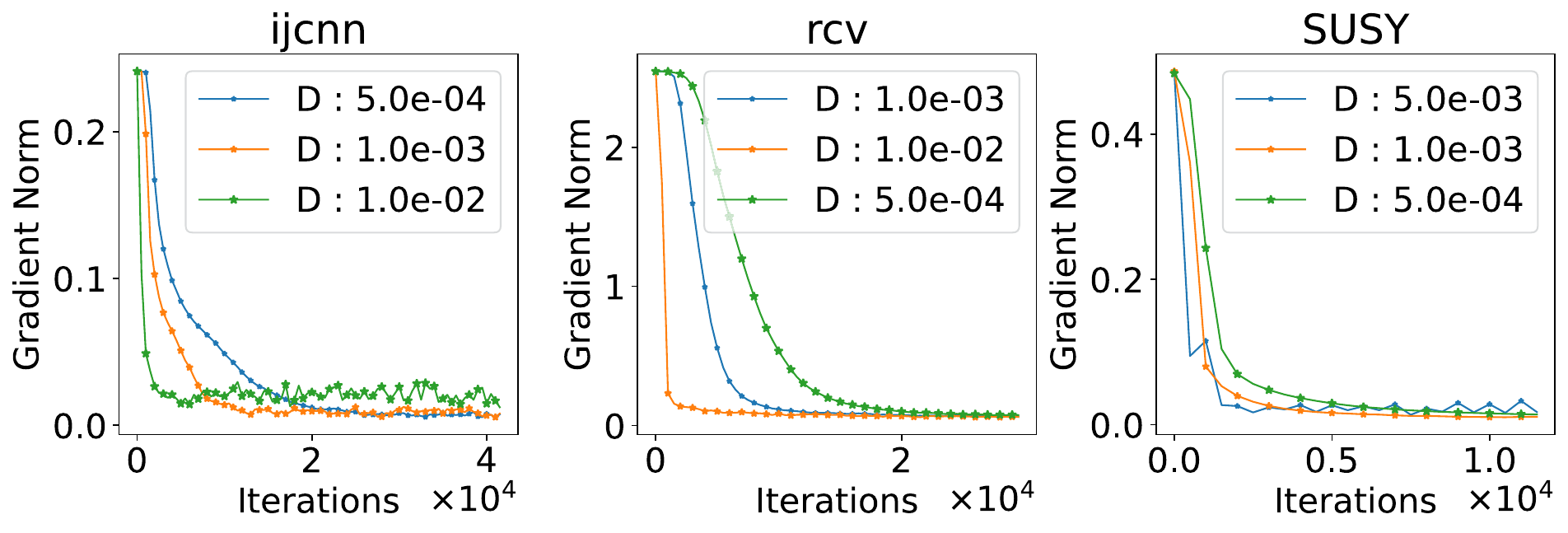}}
\caption{Evaluation of the gradient norm in the zero-order setting.}
\label{fig_grad_zo}
\end{center}
\end{figure}

{\bf{Results.}} To empirically analyze the performance of ME-DOL, global gradient norms $\norm{\nabla f(\bar{w}^k)}$ for $k\geq 1$ are calculated for both first-order and zero-order settings. We use three datasets (ijcnn, rcv, SUSY) to illustrate the decay of gradient norms with respect to iterations. The plots are reported in Figs. \ref{fig_grad_fo} and \ref{fig_grad_zo}. This observation validates Theorems \ref{thm:ns_fo} and \ref{thm:ns_zo} in our paper, respectively.

% \begin{figure}[t]
% %\vskip 0.05in

% %\vskip -0.33in
% \end{figure}

\begin{figure}[b]
%\vskip 0.05in
\begin{center}
\centerline{\includegraphics[width=\columnwidth]{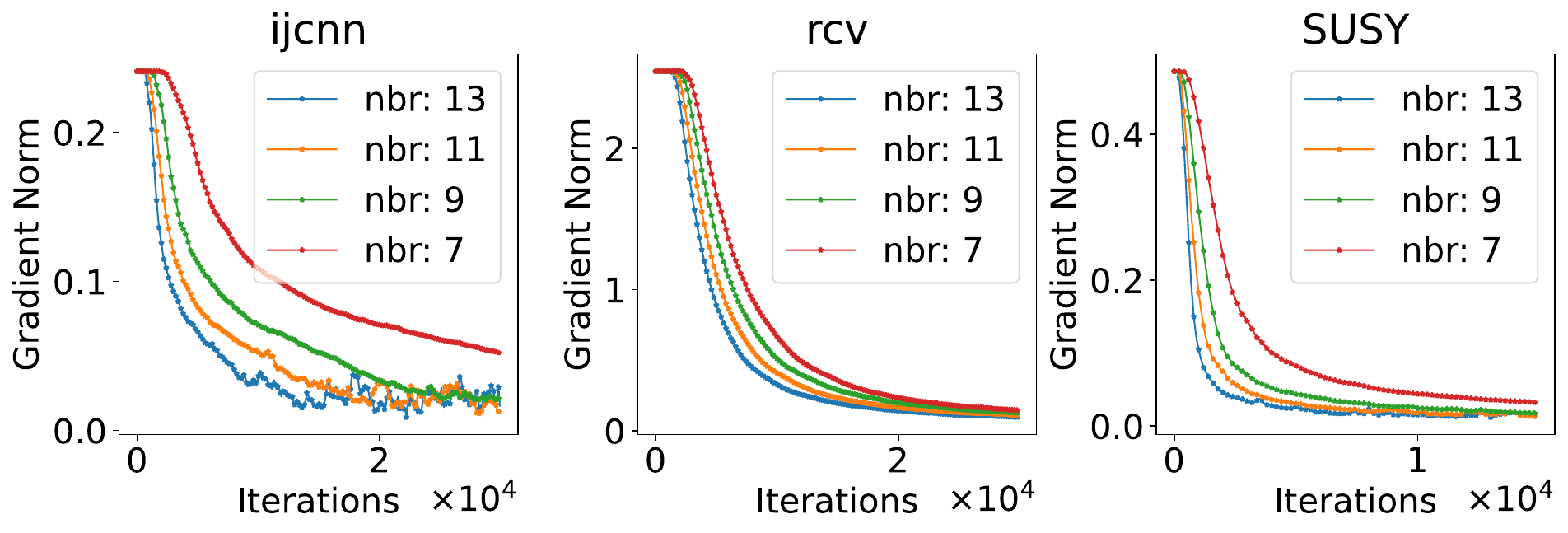}}
\caption{Ring graphs in the first-order setting.}
\label{fig_net_rho}
\end{center}
\vskip -0.18in
\end{figure}

\begin{figure}[t]
%\vskip -.43in
\begin{center}
\centerline{\includegraphics[width=\columnwidth]{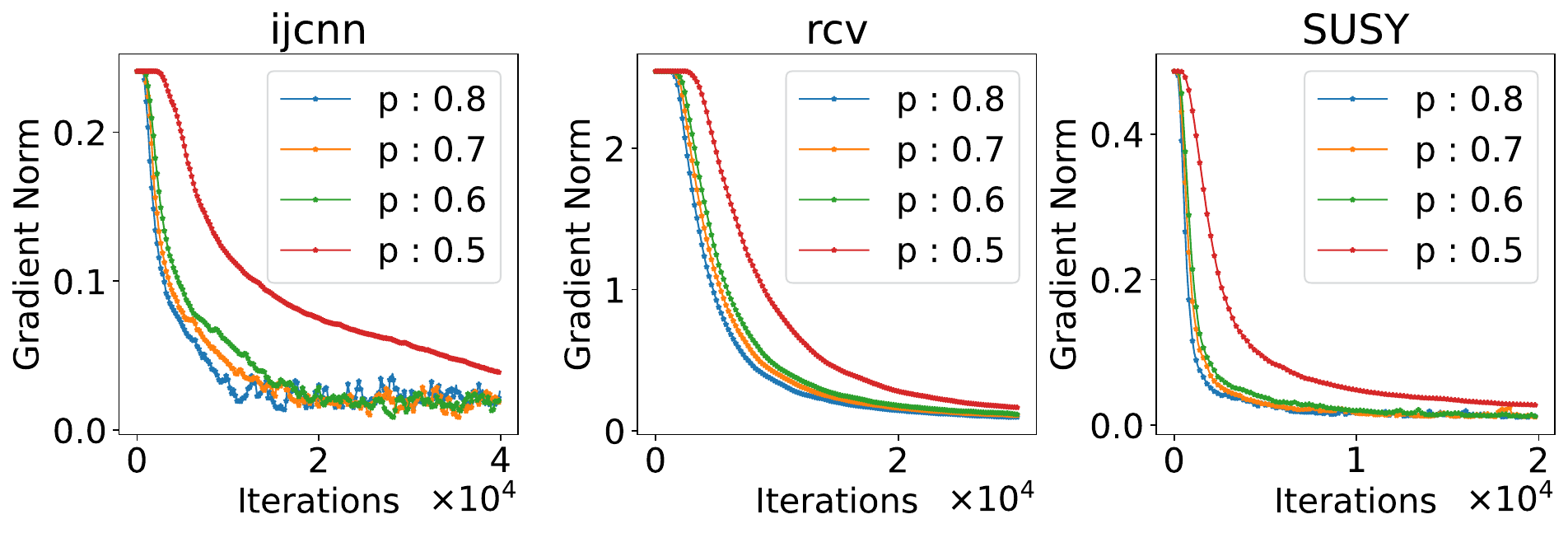}}
\caption{Random graphs in the first-order setting.}
\label{fig_netw_nosup}
\end{center}
\vskip -0.5in
\end{figure}

We further evaluate the classification accuracy over the test data. We compare our algorithm in the zero-order
setting with DGFM in \citet{lin2023decentralized} on three datasets (a9a, HIGGS, covtype),
and accuracy plots are reported in Fig. \ref{fig_comp_zero} (see Appendix \ref{app:experiment}). We can see that our algorithm
dominates DGFM in terms of the test classification accuracy.
In the first-order setting, we compare our algorithm with DPSGD in \citet{lian2017can} on three datasets (a9a,
HIGGS, covtype), and
accuracy plots are depicted in Fig. \ref{fig_comp_first} (see Appendix \ref{app:experiment}), where again our algorithm achieves a better performance in terms
of the classification accuracy. Note that DPSGD was originally proposed for smooth problems, but the same algorithm with projection (DPSM) was analyzed in the nonsmooth setting as well \cite{Chen_Distributed_weaklyconvex}.

{\bf Impact of Network:} We also evaluate the effect of network connectivity on the ring-based graphs with $n=20$ agents, 
using the number of neighbors from $\{7,9,11,13\}$. The corresponding $\rho$ values are $\{0.81, 0.70 , 0.57, 0.44 \}$, respectively. As the number of neighbors increases, the graph becomes more connected, and the value of $\rho$ decreases. In Fig. \ref{fig_net_rho} we observe that better connectivity (smaller $\rho$) results in a faster convergence in the first-order setting.

To further evaluate the effect of network topology, we design the communication matrices based on Erdos-Renyi random graph $G(20,p)$, where $p$ represents the probability of existence of an edge. In this experiment, $p$ is selected from $\{0.5,0.6,0.7,0.8\}$, for which the corresponding $\rho$ values are $\{0.83, 0.63, 0.56, 0.47\}$. The plots are presented in Fig. \ref{fig_netw_nosup}, where again we observe that larger probability, which implies (possibly) better connectivity, results in a faster convergence. Note that since in this experiment the graphs are generated randomly, one might get different values for $\rho$ even by trying the same edge probabilities.

\section{Conclusion}
We presented a novel algorithm for decentralized nonsmooth nonconvex stochastic optimization in first-order and zero-order oracle settings. We adopted recent techniques on online-to-nonconvex conversion \cite{cutkosky2023optimal} and the geometric lemma on Goldstein subdifferential sets \cite{kornowski2023algorithm} to streamline the finite-time analysis of the algorithm. Our algorithm achieved the optimal sample complexity of  $O(\delta^{-1} \epsilon^{-3})$ for finding \goldstat points of the global objective in three settings, namely (i) smooth first-order, (ii) nonsmooth first-order, and (iii) nonsmooth zero-order. Notably, to the best of our knowledge, we provided the first finite-time convergence characterization in the nonsmooth first-order setting (without weak-convexity assumption \cite{Chen_Distributed_weaklyconvex}), and our result on the nonsmooth zero-order setting does not use variance reduction. Future directions include the investigation of high probability bounds (as opposed to expectation), the optimal dependence to network parameters, as well as convergence in the deterministic regime.

In our theorems, we found that $N=O((1-\rho)^{-2})$, but it is challenging to evaluate the optimality with respect to $\rho$ in the nonsmooth setting using the ($\delta, \epsilon$)-stationarity concept. There is currently no lower bound on the communication complexity for the decentralized nonsmooth nonconvex stochastic optimization, i.e., in the nonsmooth setting the optimal dependence on $\rho$ in finding ($\delta, \epsilon$)-stationary points has not been explored yet. For the smooth decentralized setting, the lower bound on $\rho$-dependency is given as $O((1-\rho)^{-\frac{1}{2}})$ \cite{lu2021optimal}, using a carefully designed communication protocol that allows for network structure change. Whether the optimal dependence to $\rho$ in the nonsmooth setting is the same and whether that potential gap can be closed are interesting research questions.

\section*{Impact Statement}
We do not anticipate any future societal consequences as this work contributes to the theory of decentralized optimization.

\section*{Acknowledgements}
The authors gratefully acknowledge the support of Mechanical and Industrial Engineering (MIE) Chair Fellowship at Northeastern University as well as NSF ECCS-2240788 Award for this research.

\bibliography{references}

\begin{thebibliography}{38}
\providecommand{\natexlab}[1]{#1}
\providecommand{\url}[1]{\texttt{#1}}
\expandafter\ifx\csname urlstyle\endcsname\relax
  \providecommand{\doi}[1]{doi: #1}\else
  \providecommand{\doi}{doi: \begingroup \urlstyle{rm}\Url}\fi

\bibitem[Arjevani et~al.(2023)Arjevani, Carmon, Duchi, Foster, Srebro, and
  Woodworth]{arjevani2023lower}
Arjevani, Y., Carmon, Y., Duchi, J.~C., Foster, D.~J., Srebro, N., and
  Woodworth, B.
\newblock Lower bounds for non-convex stochastic optimization.
\newblock \emph{Mathematical Programming}, 199\penalty0 (1-2):\penalty0
  165--214, 2023.

\bibitem[Carmon et~al.(2020)Carmon, Duchi, Hinder, and
  Sidford]{carmon2020lower}
Carmon, Y., Duchi, J.~C., Hinder, O., and Sidford, A.
\newblock Lower bounds for finding stationary points {I}.
\newblock \emph{Mathematical Programming}, 184\penalty0 (1-2):\penalty0
  71--120, 2020.

\bibitem[Chen et~al.(2023)Chen, Xu, and Luo]{chen2023faster}
Chen, L., Xu, J., and Luo, L.
\newblock Faster gradient-free algorithms for nonsmooth nonconvex stochastic
  optimization.
\newblock In \emph{International Conference on Machine Learning (ICML)}, pp.\
  5219--5233. PMLR, 2023.

\bibitem[Chen et~al.(2021)Chen, Garcia, and
  Shahrampour]{Chen_Distributed_weaklyconvex}
Chen, S., Garcia, A., and Shahrampour, S.
\newblock On distributed nonconvex optimization: Projected subgradient method
  for weakly convex problems in networks.
\newblock \emph{IEEE Transactions on Automatic Control}, 67\penalty0
  (2):\penalty0 662--675, 2021.

\bibitem[Clarke et~al.(2008)Clarke, Ledyaev, Stern, and
  Wolenski]{clarke2008nonsmooth}
Clarke, F.~H., Ledyaev, Y.~S., Stern, R.~J., and Wolenski, P.~R.
\newblock \emph{Nonsmooth analysis and control theory}, volume 178.
\newblock Springer Science \& Business Media, 2008.

\bibitem[Cutkosky et~al.(2023)Cutkosky, Mehta, and
  Orabona]{cutkosky2023optimal}
Cutkosky, A., Mehta, H., and Orabona, F.
\newblock Optimal stochastic non-smooth non-convex optimization through
  online-to-non-convex conversion.
\newblock In \emph{International Conference on Machine Learning (ICML)}, pp.\
  6643--6670. PMLR, 2023.

\bibitem[Daniilidis \& Drusvyatskiy(2020)Daniilidis and
  Drusvyatskiy]{daniilidis2020pathological}
Daniilidis, A. and Drusvyatskiy, D.
\newblock Pathological subgradient dynamics.
\newblock \emph{SIAM Journal on Optimization}, 30\penalty0 (2):\penalty0
  1327--1338, 2020.

\bibitem[Davis \& Drusvyatskiy(2019)Davis and
  Drusvyatskiy]{davis2019stochastic}
Davis, D. and Drusvyatskiy, D.
\newblock Stochastic model-based minimization of weakly convex functions.
\newblock \emph{SIAM Journal on Optimization}, 29\penalty0 (1):\penalty0
  207--239, 2019.

\bibitem[Davis et~al.(2022)Davis, Drusvyatskiy, Lee, Padmanabhan, and
  Ye]{davis2022gradient}
Davis, D., Drusvyatskiy, D., Lee, Y.~T., Padmanabhan, S., and Ye, G.
\newblock A gradient sampling method with complexity guarantees for lipschitz
  functions in high and low dimensions.
\newblock \emph{Advances in Neural Information Processing Systems (NeurIPS)},
  35:\penalty0 6692--6703, 2022.

\bibitem[Ghadimi \& Lan(2013)Ghadimi and Lan]{ghadimi2013stochastic}
Ghadimi, S. and Lan, G.
\newblock Stochastic first-and zeroth-order methods for nonconvex stochastic
  programming.
\newblock \emph{SIAM Journal on Optimization}, 23\penalty0 (4):\penalty0
  2341--2368, 2013.

\bibitem[Jordan et~al.(2023)Jordan, Kornowski, Lin, Shamir, and
  Zampetakis]{jordan2023deterministic}
Jordan, M., Kornowski, G., Lin, T., Shamir, O., and Zampetakis, M.
\newblock Deterministic nonsmooth nonconvex optimization.
\newblock In \emph{The Thirty Sixth Annual Conference on Learning Theory}, pp.\
   4570--4597. PMLR, 2023.

\bibitem[Jordan et~al.(2022)Jordan, Lin, and Zampetakis]{jordan2022complexity}
Jordan, M.~I., Lin, T., and Zampetakis, M.
\newblock On the complexity of deterministic nonsmooth and nonconvex
  optimization.
\newblock \emph{arXiv preprint arXiv:2209.12463}, 2022.

\bibitem[Koloskova et~al.(2019)Koloskova, Stich, and Jaggi]{koloskova19a}
Koloskova, A., Stich, S., and Jaggi, M.
\newblock Decentralized stochastic optimization and gossip algorithms with
  compressed communication.
\newblock In \emph{International Conference on Machine Learning (ICML)}, pp.\
  3478--3487, 2019.

\bibitem[Kornowski \& Shamir(2021)Kornowski and Shamir]{kornowski2022oracle}
Kornowski, G. and Shamir, O.
\newblock Oracle complexity in nonsmooth nonconvex optimization.
\newblock \emph{Advances in Neural Information Processing Systems (NeurIPS)},
  34:\penalty0 324--334, 2021.

\bibitem[Kornowski \& Shamir(2024)Kornowski and Shamir]{kornowski2023algorithm}
Kornowski, G. and Shamir, O.
\newblock An algorithm with optimal dimension-dependence for zero-order
  nonsmooth nonconvex stochastic optimization.
\newblock \emph{Journal of Machine Learning Research}, 25\penalty0
  (122):\penalty0 1--14, 2024.

\bibitem[Li et~al.(2020)Li, Sahu, Zaheer, Sanjabi, Talwalkar, and
  Smith]{li2020federated}
Li, T., Sahu, A.~K., Zaheer, M., Sanjabi, M., Talwalkar, A., and Smith, V.
\newblock Federated optimization in heterogeneous networks.
\newblock \emph{Proceedings of Machine learning and systems}, 2:\penalty0
  429--450, 2020.

\bibitem[Lian et~al.(2017)Lian, Zhang, Zhang, Hsieh, Zhang, and
  Liu]{lian2017can}
Lian, X., Zhang, C., Zhang, H., Hsieh, C.-J., Zhang, W., and Liu, J.
\newblock Can decentralized algorithms outperform centralized algorithms? a
  case study for decentralized parallel stochastic gradient descent.
\newblock \emph{Advances in neural information processing systems (NeurIPS)},
  30, 2017.

\bibitem[Lin et~al.(2022)Lin, Zheng, and Jordan]{lin2022gradient}
Lin, T., Zheng, Z., and Jordan, M.
\newblock Gradient-free methods for deterministic and stochastic nonsmooth
  nonconvex optimization.
\newblock \emph{Advances in Neural Information Processing Systems (NeurIPS)},
  35:\penalty0 26160--26175, 2022.

\bibitem[Lin et~al.(2024)Lin, Xia, Deng, and Luo]{lin2023decentralized}
Lin, Z., Xia, J., Deng, Q., and Luo, L.
\newblock Decentralized gradient-free methods for stochastic non-smooth
  non-convex optimization.
\newblock In \emph{Proceedings of the AAAI Conference on Artificial
  Intelligence}, volume~38, pp.\  17477--17486, 2024.

\bibitem[Liu \& Liu(2001)Liu and Liu]{liu2001monte}
Liu, J.~S. and Liu, J.~S.
\newblock \emph{Monte Carlo strategies in scientific computing}, volume~75.
\newblock Springer, 2001.

\bibitem[Lu \& De~Sa(2021)Lu and De~Sa]{lu2021optimal}
Lu, Y. and De~Sa, C.
\newblock Optimal complexity in decentralized training.
\newblock In \emph{International Conference on Machine Learning (ICML)}, pp.\
  7111--7123. PMLR, 2021.

\bibitem[Majewski et~al.(2018)Majewski, Miasojedow, and
  Moulines]{Majewski2018AnalysisON}
Majewski, S., Miasojedow, B., and Moulines, E.
\newblock Analysis of nonsmooth stochastic approximation: the differential
  inclusion approach.
\newblock \emph{arXiv preprint arXiv:1805.01916}, 2018.

\bibitem[McMahan et~al.(2016)McMahan, Yu, Richtarik, Suresh, and
  Bacon]{mcmahan2016federated}
McMahan, H.~B., Yu, F., Richtarik, P., Suresh, A., and Bacon, D.
\newblock Federated learning: Strategies for improving communication
  efficiency.
\newblock In \emph{Proceedings of the 29th Conference on Neural Information
  Processing Systems (NeurIPS), Barcelona, Spain}, pp.\  5--10, 2016.

\bibitem[Murray et~al.(2019)Murray, Swenson, and Kar]{murray2019revisiting}
Murray, R., Swenson, B., and Kar, S.
\newblock Revisiting normalized gradient descent: Fast evasion of saddle
  points.
\newblock \emph{IEEE Transactions on Automatic Control}, 64\penalty0
  (11):\penalty0 4818--4824, 2019.

\bibitem[Nedic \& Ozdaglar(2009)Nedic and Ozdaglar]{nedic09distributed}
Nedic, A. and Ozdaglar, A.
\newblock Distributed subgradient methods for multi-agent optimization.
\newblock \emph{IEEE Transactions on Automatic Control}, 54\penalty0
  (1):\penalty0 48--61, 2009.
\newblock \doi{10.1109/TAC.2008.2009515}.

\bibitem[Rockafellar \& Wets(2009)Rockafellar and
  Wets]{rockafellar2009variational}
Rockafellar, R.~T. and Wets, R. J.-B.
\newblock \emph{Variational analysis}, volume 317.
\newblock Springer Science \& Business Media, 2009.

\bibitem[Scaman et~al.(2018)Scaman, Bach, Bubeck, Massouli{\'e}, and
  Lee]{scaman2018optimal}
Scaman, K., Bach, F., Bubeck, S., Massouli{\'e}, L., and Lee, Y.~T.
\newblock Optimal algorithms for non-smooth distributed optimization in
  networks.
\newblock \emph{Advances in Neural Information Processing Systems (NeurIPS)},
  31, 2018.

\bibitem[Scaman et~al.(2020)Scaman, Dos~Santos, Barlier, and
  Colin]{ScamanFasterSmoothing}
Scaman, K., Dos~Santos, L., Barlier, M., and Colin, I.
\newblock A simple and efficient smoothing method for faster optimization and
  local exploration.
\newblock \emph{Advances in Neural Information Processing Systems (NeurIPS)},
  33:\penalty0 6503--6513, 2020.

\bibitem[Shahrampour \& Jadbabaie(2018)Shahrampour and Jadbabaie]{ShahDecMD}
Shahrampour, S. and Jadbabaie, A.
\newblock Distributed online optimization in dynamic environments using mirror
  descent.
\newblock \emph{IEEE Transactions on Automatic Control}, 63\penalty0
  (3):\penalty0 714--725, 2018.

\bibitem[Shamir(2017)]{shamir2017optimal}
Shamir, O.
\newblock An optimal algorithm for bandit and zero-order convex optimization
  with two-point feedback.
\newblock \emph{The Journal of Machine Learning Research}, 18\penalty0
  (1):\penalty0 1703--1713, 2017.

\bibitem[Swenson et~al.(2022)Swenson, Murray, Poor, and
  Kar]{swenson2022distributed}
Swenson, B., Murray, R., Poor, H.~V., and Kar, S.
\newblock Distributed stochastic gradient descent: Nonconvexity, nonsmoothness,
  and convergence to local minima.
\newblock \emph{Journal of Machine Learning Research}, 23\penalty0
  (328):\penalty0 1--62, 2022.

\bibitem[Tang et~al.(2018)Tang, Lian, Yan, Zhang, and Liu]{tang2018d}
Tang, H., Lian, X., Yan, M., Zhang, C., and Liu, J.
\newblock $d^2$: Decentralized training over decentralized data.
\newblock In \emph{International Conference on Machine Learning (ICML)}, pp.\
  4848--4856. PMLR, 2018.

\bibitem[Tian et~al.(2022)Tian, Zhou, and So]{tian2022finitetime}
Tian, L., Zhou, K., and So, A. M.-C.
\newblock On the finite-time complexity and practical computation of
  approximate stationarity concepts of lipschitz functions.
\newblock In \emph{International Conference on Machine Learning (ICML)}, pp.\
  21360--21379. PMLR, 2022.

\bibitem[Wang et~al.(2020)Wang, Liu, Liang, Joshi, and Poor]{wang2020tackling}
Wang, J., Liu, Q., Liang, H., Joshi, G., and Poor, H.~V.
\newblock Tackling the objective inconsistency problem in heterogeneous
  federated optimization.
\newblock \emph{Advances in neural information processing systems (NeurIPS)},
  33:\penalty0 7611--7623, 2020.

\bibitem[Xu et~al.(2023)Xu, Qu, Xiang, and Gao]{xu2023asynchronous}
Xu, C., Qu, Y., Xiang, Y., and Gao, L.
\newblock Asynchronous federated learning on heterogeneous devices: A survey.
\newblock \emph{Computer Science Review}, 50:\penalty0 100595, 2023.

\bibitem[Yousefian et~al.(2012)Yousefian, Nedi{\'c}, and
  Shanbhag]{yousefian2012stochastic}
Yousefian, F., Nedi{\'c}, A., and Shanbhag, U.~V.
\newblock On stochastic gradient and subgradient methods with adaptive
  steplength sequences.
\newblock \emph{Automatica}, 48\penalty0 (1):\penalty0 56--67, 2012.

\bibitem[Zhang et~al.(2020{\natexlab{a}})Zhang, He, Sra, and
  Jadbabaie]{zhang2019gradient}
Zhang, J., He, T., Sra, S., and Jadbabaie, A.
\newblock Why gradient clipping accelerates training: A theoretical
  justification for adaptivity.
\newblock In \emph{International Conference on Learning Representations
  (ICLR)}, 2020{\natexlab{a}}.

\bibitem[Zhang et~al.(2020{\natexlab{b}})Zhang, Lin, Jegelka, Sra, and
  Jadbabaie]{zhang2020complexity}
Zhang, J., Lin, H., Jegelka, S., Sra, S., and Jadbabaie, A.
\newblock Complexity of finding stationary points of nonconvex nonsmooth
  functions.
\newblock In \emph{International Conference on Machine Learning (ICML)}, pp.\
  11173--11182, 2020{\natexlab{b}}.

\end{thebibliography}
\bibliographystyle{icml2024}

%%%%%%%%%%%%%%%%%%%%%%%%%%%%%%%%%%%%%%%%%%%%%%%%%%%%%%%%%%%%%%%%%%%%%%%%%%%%%%%
%%%%%%%%%%%%%%%%%%%%%%%%%%%%%%%%%%%%%%%%%%%%%%%%%%%%%%%%%%%%%%%%%%%%%%%%%%%%%%%
% APPENDIX
%%%%%%%%%%%%%%%%%%%%%%%%%%%%%%%%%%%%%%%%%%%%%%%%%%%%%%%%%%%%%%%%%%%%%%%%%%%%%%%
%%%%%%%%%%%%%%%%%%%%%%%%%%%%%%%%%%%%%%%%%%%%%%%%%%%%%%%%%%%%%%%%%%%%%%%%%%%%%%%
\newpage
\appendix
\onecolumn
%\section{You \emph{can} have an appendix here.}

%You can have as much text here as you want. The main body must be at most $8$ pages long. For the final version, one more page can be added. If you want, you can use an appendix like this one. The $\mathtt{\backslash onecolumn}$ command above can be kept in place if you prefer a one-column appendix, or can be removed if you prefer a two-column appendix.  Apart from this possible change, the style (font size, spacing, margins, page numbering, etc.) should be kept the same as the main body.
%%%%%%%%%%%%%%%%%%%%%%%%%%%%%%%%%%%%%%%%%%%%%%%%%%%%%%%%%%%%%%%%%%%%%%%%%%%%%%%
%%%%%%%%%%%%%%%%%%%%%%%%%%%%%%%%%%%%%%%%%%%%%%%%%%%%%%%%%%%%%%%%%%%%%%%%%%%%%%%

\section{Proof of Theorems }
First, we state the following lemma to use in the proof of Theorem \ref{thm:bigtheorem}.
\begin{lemma}\label{lem:dec_bound}
Given Assumption \ref{assumption:dsm}, for the update in Algorithm \ref{alg:decen}, we have for any $i\in [n]$, $t\in [T]$, $k\in [K]$ that
$$
\norm{\bar{y}^k_t-y^k_{t,i}}\leq \frac{D\sqrt{n}}{1-\rho}.
$$
\begin{proof}
We have based on the update rule that
\begin{align*}
y^k_{t,i} = \sum_{j=1}^nP_{ij}x^k_{t,j} = \sum_{j=1}^nP_{ij} y_{t-1,j}^k+\sum_{j=1}^nP_{ij}\Delta_{t,j}^k.
\end{align*}
Let $y^k_{t} \in \mathbb{R}^{nd}$ be the concatenation of vectors $y^k_{t,1},y^k_{t,2},\ldots,y^k_{t,n}$, and $\zeta^k_{t} \in \mathbb{R}^{nd}$ be the concatenation of vectors $\sum_{j=1}^nP_{1j}\Delta^k_{t,j},\sum_{j=1}^nP_{2j}\Delta^k_{t,j},\ldots,\sum_{j=1}^nP_{nj}\Delta^k_{t,j}$. We then have
$$
y^k_t=(P\otimes I) y^k_{t-1} + \zeta_t^k.
$$
Without loss of generality let $y^k_{0}=0$. Then,
$$
y^k_{t,i}=\sum_{j=1}^n\sum_{\tau=0}^{t-1}[P^{t-1-\tau}]_{ij}\zeta^k_{\tau+1,j} \Rightarrow  y^k_{t,i}-\bar{y}^k_t=\sum_{j=1}^n\sum_{\tau=0}^{t-1}\Big([P^{t-1-\tau}]_{ij}-\frac{1}{n}\Big)\zeta^k_{\tau+1,j}.
$$
Combining the geometric mixing bound of $\sum_{j=1}^{n}|P_{ij}^t-\frac{1}{n}| \le \sqrt{n} \rho^t$ \cite{liu2001monte} and the fact that $\|\zeta^k_{\tau+1,j}\|\le D$, the proof is complete. 
\end{proof}
\end{lemma}

\subsection{Proof of Proposition \ref{prop:deltanorm}}
\begin{proof}
\label{proof:propdeltanorm}

By Lemma \ref{lem:geo_gold} we have $\partial_{\mu} f_{\delta}(x) \subseteq \partial_{\mu+\delta}f(x)$. Using  Definition \ref{def:goldstein}, we have  $\norm{\nabla f(x)}_{\mu +\delta} \le \norm{\nabla f_{\delta}(x)}_{\mu} $. Replacing $\delta$ with $a\delta$ and $\mu$ with $(1-a)\delta$ for $a\in (0,1)$ gives $\norm{\nabla f(x)}_{\delta} \le \norm{\nabla f_{a\delta}(x)}_{(1-a)\delta} $.
\end{proof}

\subsection{Proof of Lemma \ref{lem:decentralized}}
\begin{proof}
We know that $$\frac{1}{n} \sum_{i=1}^{n}\nabla f^{i}(\bar{w}_t) = \nabla f(\bar{w}_t) = \frac{1}{n} \sum_{i=1}^{n}\nabla f(\bar{w}_t).$$
Using $L_1$ smoothness of $f^i$ and $f$ and $\norm{w_{t,i} -\bar{w}_t} \le r$, we have that $\norm{\nabla f^{i}(w_{t,i}) - \nabla f^{i}(\bar{w}_t) } \le r L_1$ and $\norm{\nabla f(w_{t,i}) - \nabla f(\bar{w}_t) }\le r L_1$. Therefore, 
\begin{align*}
\norm{\frac{1}{nT} \finitesum{t}{T} \finitesum{i}{n} \nabla f(w_{t,i}) } 
% & = 
%  \norm{ \frac{1}{nT} \finitesum{t}{T} \left( \sum_{i=1}^{n}\nabla f(w_{t,i})
% -  \sum_{i=1}^{n}\nabla f(\bar{w}_t)
% +  \sum_{i=1}^{n}\nabla f^{i}(\bar{w}_t)
% -  \sum_{i=1}^{n}\nabla f^{i}(w_{t,i}) 
% +  \sum_{i=1}^{n}\nabla f^{i}(w_{t,i}) \right) } \\
     & \le
     \frac{1}{nT} \finitesum{t}{T} \norm{ \sum_{i=1}^{n}\nabla f(w_{t,i})
    - \sum_{i=1}^{n}\nabla f(\bar{w}_t)}
    + \frac{1}{nT} \finitesum{t}{T} \norm{ \sum_{i=1}^{n}\nabla f^{i}(\bar{w}_t)
    - \sum_{i=1}^{n}\nabla f^{i}(w_{t,i}) } \\
    & + \norm{\frac{1}{nT} \finitesum{t}{T}  \sum_{i=1}^{n}\nabla f^{i}(w_{t,i}) }
    \\
     & \le \norm{\frac{1}{nT} \finitesum{t}{T} \finitesum{i}{n} \nabla f^i(w_{t,i}) } +2rL_1,
\end{align*}
which completes the proof.
\end{proof}

\subsection{Proof of Theorem \ref{thm:bigtheorem}}\label{proof:Thm31}
\begin{proof}
Throughout the proof, superscript $k$ denotes the $k$-th epoch, subscript $i$ denotes the agent index, and subscript $t$ represents the iteration. 

Let us start with the following definitions:

$g_{t,i}^k := \mathcal{O}_{f}^i(w_{t,i}^k)$ (stochastic gradient returned by the oracle)

$ \bar{g}_t^k := \frac{1}{n}\sum_{i=1}^{n} g_{t,i}^k $ (average of local stochastic gradients)

$\Delta_{t,i}^k$ is generated based on the decentralized online learning algorithm $\mathcal{A}$ (Algorithm \ref{alg:DOL})

$ \bar{\Delta}_t^k := \frac{1}{n}\sum_{i=1}^{n} \Delta_{t,i}^k $ (average of local actions)

$\nabla_{t,i}^k := \int_{0}^{1}\nabla f^{i}(y_{t-1,i}^k+s\Delta_{t,i}^k)ds$ (expected local gradient)

$ \bar{\nabla}_t^k := \frac{1}{n}\sum_{i=1}^{n} \nabla_{t,i}^k $ (average of expected local gradients)

$\tilde{\nabla}_t^k :=  \int_{0}^{1}\nabla f(\bar{x}_{t-1}^k+s\bar{\Delta}_{t}^k)ds $  (expected global gradient)

For any Lipschitz continuous function and an update rule $x_t = x_{t-1} + \Delta_t$ we have
$$ f(x_t) - f(x_{t-1}) = \int_{0}^{1} \innerproduct{\nabla f(x_{t-1} + s\Delta_{t}),  \Delta_t} ds = \innerproduct{\nabla_t, \Delta_t}.$$
In our decentralized update rule it follows by doubly stochasticity of $P$ that 
$$\bar{y}_t^k:=\frac{1}{n}\sum_{i=1}^n y_{t,i}^k = \frac{1}{n}\sum_{i=1}^n\sum_{j=1}^n P_{ij} x_{t,j}^k= \frac{1}{n}\sum_{j=1}^n  \sum_{i=1}^n P_{ij} x_{t,j}^k= \frac{1}{n}\sum_{j=1}^n x_{t,j}^k=:\bar{x}_t^k,$$
and since $x_{t,i}^k = y_{t-1,i}^k+\Delta_{t,i}^k$,  we get that $\bar{x}_t^k = \bar{x}^k_{t-1} + \bar{\Delta}^k_t$.

Therefore, we can write the following for the global function $f(x) = \frac{1}{n} \sum_{i=1}^{n}f^i(x)$,
$$f(\bar{x}_t^k) - f(\bar{x}_{t-1}^k) =  \innerproduct{ \tilde{\nabla}_t^k ,  \bar{\Delta}_t^k},$$
and summing both sides over $t\in[T]$ gives $$f(\bar{x}_T^k) - f(\bar{x}_{0}^k) =  \finitesum{i}{T} \innerproduct{ \tilde{\nabla}_t^k , \bar{\Delta}_t^k }.$$
There are two sources of randomness in the algorithm, namely $\xi_{t,i}^k$ and $s_{t,i}^k$. Taking expectation over those, we can decompose above into four terms: 
\begin{align}\label{eq:first}
\mathbb{E}[f(\bar{x}_T^k)- f(\bar{x}_0^k)] &= \sum_{t=1}^{T} \mathbb{E}[\innerproduct{\bar{\Delta}_{t}^k ,  \tilde{\nabla}_t^k}] \notag\\
&=  \underbrace{\sum_{t=1}^{T} \mathbb{E}[\innerproduct{ \bar{g}_{t}^k, \bar{\Delta}_t^k - u^k  }]}_{R_T^k(u^k)} + \underbrace{\sum_{t=1}^{T} \mathbb{E}[\innerproduct{\bar{g}_t^k, u^k}]}_{T_2} 
+ \underbrace{\sum_{t=1}^{T} \mathbb{E}[ \innerproduct{\bar{\Delta}_t^k , \tilde{\nabla}_t^k - \bar{\nabla}_t^k }]}_{T_3} 
+ \sum_{t=1}^{T} \mathbb{E}[ \innerproduct{\bar{\Delta}_t^k , \bar{\nabla}_t^k - \bar{g}_t^k }], 
\end{align}
where the last term equals zero due to the unbiased gradient assumption that $\mathbb{E}[\bar{g}_t^k]=\bar{\nabla}_t^k$. The above holds for any $u^k$, and choosing $ u^k = -D\frac{\lhslocal }{\norm{\lhslocal}}$, we have that

\begin{align*}
T_2 & =\probexp{\innerproduct{\sum_{t=1}^{T} \bar{g}_t^k, u^k}} = \probexp{\innerproduct{ u^k , \frac{1}{n}\lhslocal }} + \probexp{ \innerproduct{ u^k , \finitesum{t}{T}\bar{g}_t^k - \frac{1}{n} \lhslocal }}&\\
     & \le \probexpbig{-DT \norm{\frac{1}{nT}\lhslocal} } + \probexpbig{\frac{D}{n} \norm{\finitesum{t}{T} \finitesum{i}{n} (\nabla f^i (w_{t,i}^k) - g_{t,i}^k ) }   } &\\ 
     & \le \probexpbig{-DT \norm{\frac{1}{nT}\lhslocal} } + D\sigma\sqrt{\frac{T}{n}}, 
\end{align*}
where we used Assumption \ref{assumption:gradientoracle} and Jensen's inequality in the last line. 

Rearranging Equation \eqref{eq:first} using above and dividing by $DT$ yields 
\begin{equation}\label{eq:sec}
\underbrace{ \probexpbig{\norm{\frac{1}{nT} \lhslocal } }  }_{\epsilon-term} \le \underbrace{\frac{\mathbb{E}[f(\bar{x}^k_0) - f(\bar{x}^k_T)]}{DT}} _{sub-optimality}+\underbrace{\frac{\sigma}{\sqrt{nT}}}_{noise} + 
\underbrace{\regretterm}_{regret-term} + 
\underbrace{\discrepancyterm}_{discrepancy}. 
\end{equation}
We will now average above over $K$ epochs and bound each term. It is only the sub-optimality term that will telescope. Other terms can be bounded independent of $k$, i.e., 
bounds for noise term, regret term and discrepancy term are independent of epochs. Recall that $N:=KT$. 

\textbf{Sub-optimality Term:}
Let $\gamma := f(\bar{x}_0)-\inf_x f(x)$. Summing the sub-optimality term over $k\in [K]$ and dividing by $K$, the sum telescopes as follows due to the initialization at the start of each epoch $k\in[K]$: 
\begin{equation}\label{eq:bound1}
\frac{1}{K}\sum_{k=1}^K\frac{\mathbb{E}[f(\bar{x}^k_0) - f(\bar{x}^k_T)]}{DT}=\frac{f(\bar{x}^1_0) - \mathbb{E}[f(\bar{x}^K_T)]}{DTK} \le \frac{\gamma}{DN}.    
\end{equation}
\textbf{Regret Term:}
For the regret term, we can use Theorem 5 of \citet{ShahDecMD} with a fixed learning rate $\eta$, where for any $k\in [K]$:
$$R_T^k(u^k) \le \frac{4D^2}{\eta} + \eta \Big( \frac{G^2 T}{2} + \frac{2TG(L+G) \sqrt{n} }{1-\rho}\Big).$$ 

Choosing $\eta = \frac{8D}{c_1\sqrt{T}} $ where $c_1 = 4\sqrt{\frac{G^2(1-\rho)+4G(L+G)\sqrt{n}}{2(1-\rho)}}$ gives 
\begin{align}\label{eq:bound2}
regret-term \le \frac{D\sqrt{T}c_1}{DT} = O(T^{-1/2}).
\end{align}

\textbf{Discrepancy Term:} For this part, we remove the superscript $k$ for simplicity as the results hold for any $k\in [K]$. 
First, recall that $\|\bar{\Delta}_t\|\le D$ since the domain $\mathcal{D}$ in Algorithm \ref{alg:DOL} is bounded. Next, we will bound $\|\tilde{\nabla}_t - \bar{\nabla}_t\|$  under the assumption that local functions $f^i$ are $L_1$-Lipschitz smooth. Note that due to doubly stochasticity of $P$ we also have $\bar{x}_t = \bar{y}_t$. Therefore, 
$$\tilde{\nabla}_{t} - \bar{\nabla}_t = \frac{1}{n} \sum_{i=1}^{n} \int_{0}^{1} (\nabla f^i(\bar{y}_{t-1} + s \bar{\Delta}_{t})  - \nabla f^i(y_{t-1,i} + s \Delta_{t,i}))ds.$$
Then, we have
$$ \|\tilde{\nabla}_{t} - \bar{\nabla}_t \| \le \frac{L_1}{n} \sum_{i=1}^{n} \int_{0}^{1}\norm{\bar{y}_{t-1} + s \bar{\Delta}_{t}  - y_{t-1,i} - s \Delta_{t,i}}ds \le \frac{L_1}{n} \sum_{i=1}^{n} \norm{ \bar{y}_{t-1} - y_{t-1,i} } + \frac{L_1}{2n} \sum_{i=1}^{n} \norm{ \bar{\Delta}_{t} - \Delta_{t,i} }.$$
The first term can be bounded with Lemma \ref{lem:dec_bound} as $\norm{ \bar{y}_{t-1} - y_{t-1,i} } \le \frac{D\sqrt{n} }{1-\rho}$ . We can bound the second term with $\norm{\bar{\Delta}_{t} - \Delta_{t,i} } \le 2D$. Hence, $\|\tilde{\nabla}_{t} - \bar{\nabla}_t \| \le L_1 Dc_2$ where $c_2 := \frac{\sqrt{n}}{1-\rho} + 1$. 
%We also use the bound on   $\norm{ \bar{x}_{t-1} - x_{t-1,i} } $ later.
The discrepancy term can then be bounded as 
\begin{align}\label{eq:bound3} 
%\frac{\sum_{t=1}^{T} \innerproduct{\bar{\Delta}_t, \tilde{\nabla}_t - \bar{\nabla}_t  }}{DT}
discrepancy-term \le DL_1c_2.
\end{align}

Substituting \eqref{eq:bound1}, \eqref{eq:bound2}, and \eqref{eq:bound3} into \eqref{eq:sec}, we get
\begin{align}\label{eq:bound6}
    \epsterm \le \frac{\gamma}{DN} + \frac{\sigma}{\sqrt{nT}} + \frac{c_1}{\sqrt{T}} + DL_1c_2.
\end{align}
In the left-hand side of \eqref{eq:bound6}, we have the average of local gradients. Using Lemma \ref{lem:decentralized}, we can connect this to the average of global gradients. To this end, we need $r$ such that $\|w_{t,i}^k - \bar{w}_t^k\| \le r $. Since $\| \bar{y}_{t-1}^k - y_{t-1,i}^k\| \le \frac{D\sqrt{n} }{1-\rho}$, we have $\|w^k_{t,i}-\bar{w}_t^k\| \le 
\| \bar{y}_{t-1}^k - y_{t-1,i}^k \| + 2D \le D (\frac{\sqrt{n} }{1-\rho} +2)$. Now, utilizing Lemma \ref{lem:decentralized} with $r =D(\frac{\sqrt{n} }{1-\rho} +2) $, we obtain
\begin{align}\label{eq:bound7}
\norm{\frac{1}{nT} \finitesum{t}{T} \finitesum{i}{n} \nabla f(w_{t,i}^k) } \le \norm{\frac{1}{nT} \finitesum{t}{T} \finitesum{i}{n} \nabla f^i(w_{t,i}^k) } +DL_1 \Big(\frac{2\sqrt{n} }{1-\rho} +4\Big).
\end{align}
Combining \eqref{eq:bound6} and \eqref{eq:bound7}, we obtain $$\epstermglobal \le \frac{\gamma}{DN} + \frac{\sigma}{\sqrt{nT}} + \frac{c_1}{\sqrt{T}} + DL_1c_3,$$
where $c_3 := c_2 + \frac{2\sqrt{n} }{1-\rho} +4 = \frac{3\sqrt{n} }{1-\rho} +5$.

We must have $\norm{w_{t,i}^k -\bar{w}^k} \le \delta$ to bound $\norm{\nabla f(\bar{w}^k)}_{\delta}$. We have $\norm{w_{t_1,i}^k-w_{t_2,j}^k} \le \norm{w_{t_1,i}^k -\bar{w}_{t_1}^k} + \norm{\bar{w}_{t_1}^k - \bar{w}_{t_2}^k} + \norm{w_{t_2,j}^k -\bar{w}_{t_2}^k} \le 2r + DT = D(\frac{2\sqrt{n} }{1-\rho} +4 + T )$. If $T \ge 3 $ and $\frac{\sqrt{n}}{1-\rho}\ge 2$, choosing $D = \frac{\delta  (1-\rho)}{2T\sqrt{n}}$
guarantees that $\norm{w_{t,i}^k -\bar{w}^k} \le \delta$, and thus
\begin{align}
\mathbb{E}_{k\sim \un[K]}\left[ \gradnormdelta{f}{\bar{w}^k}\right] \le \epstermglobal \le \frac{2\gamma T \sqrt{n}  }{\delta N (1-\rho)} + \frac{\sigma}{\sqrt{nT}} + \frac{c_1}{\sqrt{T}} + \frac{\delta L_1 (1-\rho) c_3}{2T \sqrt{n}}.
\label{eq:smooth_ineq}
\end{align}
This inequality will also be used in the proof of Theorem \ref{thm:ns_fo} and Theorem \ref{thm:ns_zo}. Now, we can use $\delta < 1$ and $\frac{1}{T} \le \frac{1}{\sqrt{T}}$ to get
$$ \epstermglobal \le \frac{2\gamma \sqrt{n} }{\delta N (1-\rho)} T + \frac{1}{\sqrt{T}} \left(\frac{\sigma}{\sqrt{n}} + c_1 + L_1 \frac{1-\rho}{2\sqrt{n}} c_3\right).$$
Choosing $T = c_4 \left( \delta N \right)^{\frac{2}{3}} $ where $c_4 := \left( \frac{(1-\rho)(2 \sigma + 2c_1 \sqrt{n} + L_1 (1-\rho) c_3)}{8\gamma n} \right)^{\frac{2}{3}}$, we have
$$\epstermglobal \le (\delta N)^{-\frac{1}{3}} \left(\frac{2\gamma \sqrt{n} c_4}{ 1-\rho}  + \frac{1}{\sqrt{c_4}} \left(\frac{\sigma}{\sqrt{n}} + c_1 + L_1 \frac{1-\rho}{2\sqrt{n}} c_3\right)\right ) = c_5 (\delta N)^{-\frac{1}{3}},$$
where $c_5 := \frac{6\gamma \sqrt{n}}{1-\rho} c_4  = \frac{6 \gamma \sqrt{n}}{1-\rho} \left( \frac{(1-\rho)(2 \sigma + 2c_1 \sqrt{n} + L_1 (1-\rho) c_3)}{8\gamma n} \right)^{\frac{2}{3}} = \frac{3}{2} \left( \frac{\gamma(2 \sigma + 2c_1 \sqrt{n} + L_1 (1-\rho) c_3)^2}{(1-\rho) \sqrt{n}} \right)^{\frac{1}{3}}.$ In terms of the network connectivity measure $1-\rho$,   $c_5=O((1-\rho)^{-\frac{2}{3}})$.
As a result, we derive 
$$\mathbb{E}_{k\sim \un[K] }\left[ \gradnormdelta{f}{\bar{w}^k}\right]\le \epstermglobal \le O((\delta N)^{-\frac{1}{3}}),$$
which means that given $\{\delta,\epsilon,\rho\}$, we can find a \goldstat point in $N=\Theta\big(\delta^{-1}\epsilon^{-3}(1-\rho)^{-2}\big)$ rounds. For smooth functions $\delta = O(\epsilon)$, so the overall rate matches the optimal rate of $N=O(\epsilon^{-4})$. 
\end{proof}

\subsection{Proof of Theorem \ref{thm:ns_fo}}

\begin{proof}
\label{proof:ns_fo}
Recall from Proposition \ref{prop:ran_smo} that $(f^i)_\delta$ and in turn $f_\delta$ have $L_1$-Lipschitz smooth gradients with smoothness parameter $L_1 = \geoconstant L\sqrt{d}\delta^{-1}$, where $L$ is due to Lipschitz continuity of the original functions $f^i$ and $f$. Now, in Equation \eqref{eq:smooth_ineq} replacing $\delta$ by $\frac{\delta}{2}$ and $f$ by $f_{\frac{\delta}{2}}$, the smoothness parameter becomes $L_1 = \frac{2\geoconstant L\sqrt{d}}{\delta}$, which yields
\begin{align}
\epstermglobalf{f_{\frac{\delta}{2}}} \le \frac{4 \gamma' T \sqrt{n} }{\delta N (1-\rho)} + \frac{G}{\sqrt{nT}} + \frac{c_1}{\sqrt{T}} + \frac{\geoconstant L\sqrt{d} (1-\rho) c_3}{2T\sqrt{n}},
\label{eq:ns_fo_ineq}
\end{align}
where $\gamma':=\gamma+L$ due to the approximation error incurred in \eqref{eq:bound1}, and $\sigma$ is also replaced by $G$. For the last term we can use $\frac{1}{T} \le \frac{1}{\sqrt{T}}$ to get
$$ \epstermglobalf{f_{\frac{\delta}{2}}} \le \frac{4\gamma' \sqrt{n} }{\delta N (1-\rho)} T + \frac{1}{\sqrt{T}} \left(\frac{G}{\sqrt{n}} + c_1 + \frac{\geoconstant L\sqrt{d} (1-\rho) c_3}{2\sqrt{n}}\right).$$
Choosing $T = c_7 \left( \delta N \right)^{\frac{2}{3}} $ where $c_7 := \left( \frac{(1-\rho)(2G + 2c_1\sqrt{n} + \geoconstant L\sqrt{d} (1-\rho) c_3)}{16\gamma' n} \right)^{\frac{2}{3}}$, we have
$$ \epstermglobalf{f_{\frac{\delta}{2}}} \le (\delta N)^{-\frac{1}{3}} \left(\frac{4\gamma' \sqrt{n} c_7}{ (1-\rho)}  + \frac{1}{\sqrt{c_7}} \left(\frac{G}{\sqrt{n}} + c_1 + \frac{\geoconstant L\sqrt{d} (1-\rho)}{2\sqrt{n}} c_3\right)\right ) = c_8 (\delta N)^{-\frac{1}{3}},$$
where $c_8 :=  \frac{12\gamma' \sqrt{n}}{1-\rho}c_7  = \frac{12\gamma' \sqrt{n}}{1-\rho} \left( \frac{(1-\rho)(2G + 2c_1\sqrt{n} + \geoconstant L\sqrt{d} (1-\rho) c_3)}{16\gamma' n} \right)^{\frac{2}{3}} = 3 \left( \frac{\gamma' (2G + 2c_1\sqrt{n} + \geoconstant L\sqrt{d} (1-\rho) c_3)^2}{4(1-\rho) \sqrt{n}} \right)^{\frac{1}{3}}.$ Using Proposition \ref{prop:deltanorm} on the left-hand side of above completes the proof. In terms of the network connectivity measure $1-\rho$, we have $c_8=O((1-\rho)^{-\frac{2}{3}})$.
\end{proof}

\subsection{Proof of Theorem \ref{thm:ns_zo}}

\begin{proof}
\label{proof:ns_zo}
Similar to the proof of Theorem \ref{thm:ns_fo},  in Equation \eqref{eq:smooth_ineq}, we replace $\delta$ by $\frac{\delta}{2}$ and $f$ by $f_{\frac{\delta}{2}}$, so the smoothness parameter becomes $L_1 = \frac{2\geoconstant L\sqrt{d}}{\delta}$. Applying Lemma \ref{lemma:zeroorder}, we can bound $\sigma$ by $\sqrt{\zeroordersecbound}$ as well. Therefore, we have
$$\epstermglobalf{f_{\frac{\delta}{2}}} \le \frac{4\gamma' T \sqrt{n} }{\delta N (1-\rho)} + \frac{\sqrt{\zeroordersecbound}}{\sqrt{nT}} + \frac{c_{10}}{\sqrt{T}} + \frac{\geoconstant L\sqrt{d} (1-\rho) c_3}{2T \sqrt{n}},$$
where $c_{10}:=  4\sqrt{\frac{c_9^2(1-\rho)+4c_9(L+c_9)\sqrt{n}}{2(1-\rho)}}  $ and $c_9:= \sqrt{\zeroordersecbound}.$

If we follow similar steps as in the proof of Theorem \ref{thm:ns_fo}, we have the following result
$$ \epstermglobalf{f_{\frac{\delta}{2}}} \le  c_{11} (\delta N)^{-\frac{1}{3}},$$
where $c_{11} :=  3 \left( \frac{\gamma' (2\sqrt{\zeroordersecbound} + 2c_{10}\sqrt{n} + \geoconstant L\sqrt{d} (1-\rho) c_3)^2}{4 (1-\rho) \sqrt{n}} \right)^{\frac{1}{3}}.$ Using Proposition \ref{prop:deltanorm} on the left-hand side of above completes the proof. In terms of the network connectivity measure $1-\rho$ and ambient dimension $d$, we have    $c_{11}=O(d^{\frac{1}{3}}(1-\rho)^{-\frac{2}{3}})$.
\end{proof}

\section{Constant Terms} \label{app:constant}
The smoothness parameter of $f_\delta$ is $\kappa \frac{d!!}{(d-1)!!}\frac{L}{\delta}$, where $\kappa = \frac{2}{\pi}$ if $d$ is even, and $\kappa=1$ otherwise. Thus, the geometric constant $c := \kappa \frac{1}{\sqrt{d}}\frac{d!!}{(d-1)!!}$. We note that $\underset{d \to \infty}{\text{lim}} c = \underset{d \to \infty}{\text{lim}} \kappa \frac{1}{\sqrt{d}}\frac{d!!}{(d-1)!!} = \sqrt{\frac{\pi}{2}} $ \cite{yousefian2012stochastic}. Here, we summarize the constant terms used throughout the proofs:
\begin{align*}
c_1 &= 4\sqrt{\frac{G^2(1-\rho)+4G(L+G)\sqrt{n}}{2(1-\rho)}} \\
c_3 &= \frac{3\sqrt{n} }{1-\rho} +5\\
c_5 &= \frac{3}{2} \left( \frac{\gamma(2 \sigma + 2c_1 \sqrt{n} + L_1 (1-\rho) c_3)^2}{(1-\rho) \sqrt{n}} \right)^{\frac{1}{3}}\\
c_8 &= 3 \left( \frac{\gamma' (2G + 2c_1\sqrt{n} + \geoconstant L\sqrt{d} (1-\rho) c_3)^2}{4(1-\rho) \sqrt{n}} \right)^{\frac{1}{3}}\\
c_9 &= \sqrt{\zeroordersecbound}\\
c_{10} &=  4\sqrt{\frac{c_9^2(1-\rho)+4c_9(L+c_9)\sqrt{n}}{2(1-\rho)}} \\
c_{11} &= 3 \left( \frac{\gamma' (2\sqrt{\zeroordersecbound} + 2c_{10}\sqrt{n} + \geoconstant L\sqrt{d} (1-\rho) c_3)^2}{4(1-\rho) \sqrt{n}} \right)^{\frac{1}{3}}
\end{align*}

\section{Numerical Experiments Results}
\label{app:experiment}

\begin{figure}[t!]
\begin{center}
\vskip -3.4in
\centerline{\includegraphics[width=\columnwidth]{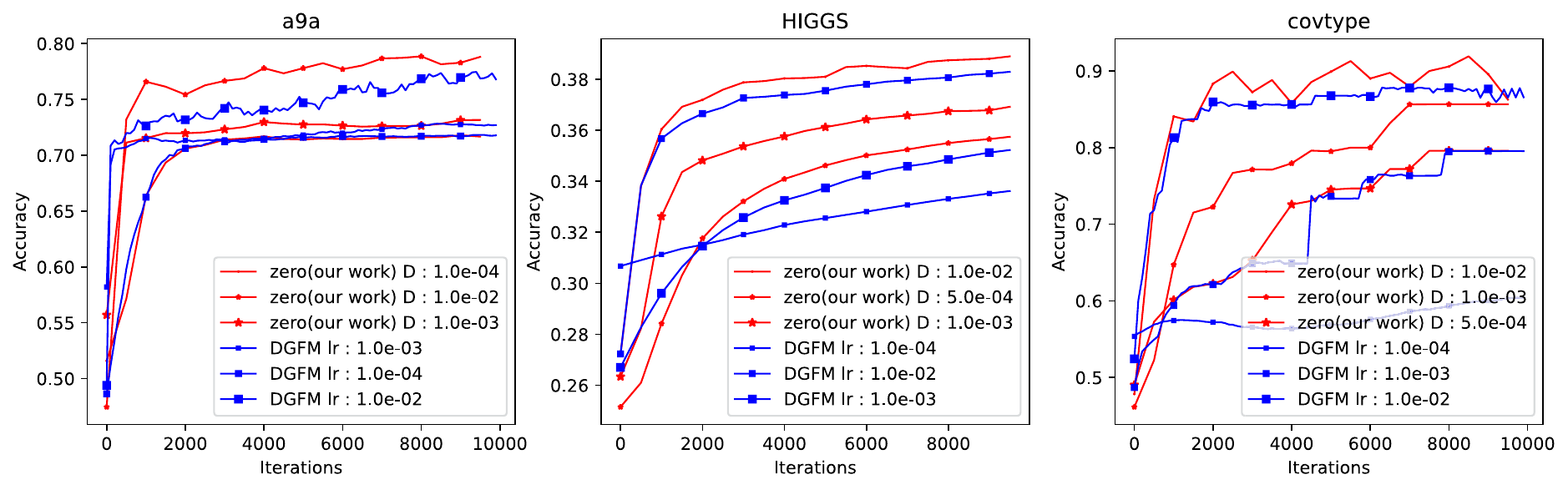}}
\caption{Evaluation of the test accuracy of our algorithm and DGFM in the zero-order setting.  }
\label{fig_comp_zero}
\end{center}
\begin{center}
\centerline{\includegraphics[width=\columnwidth]{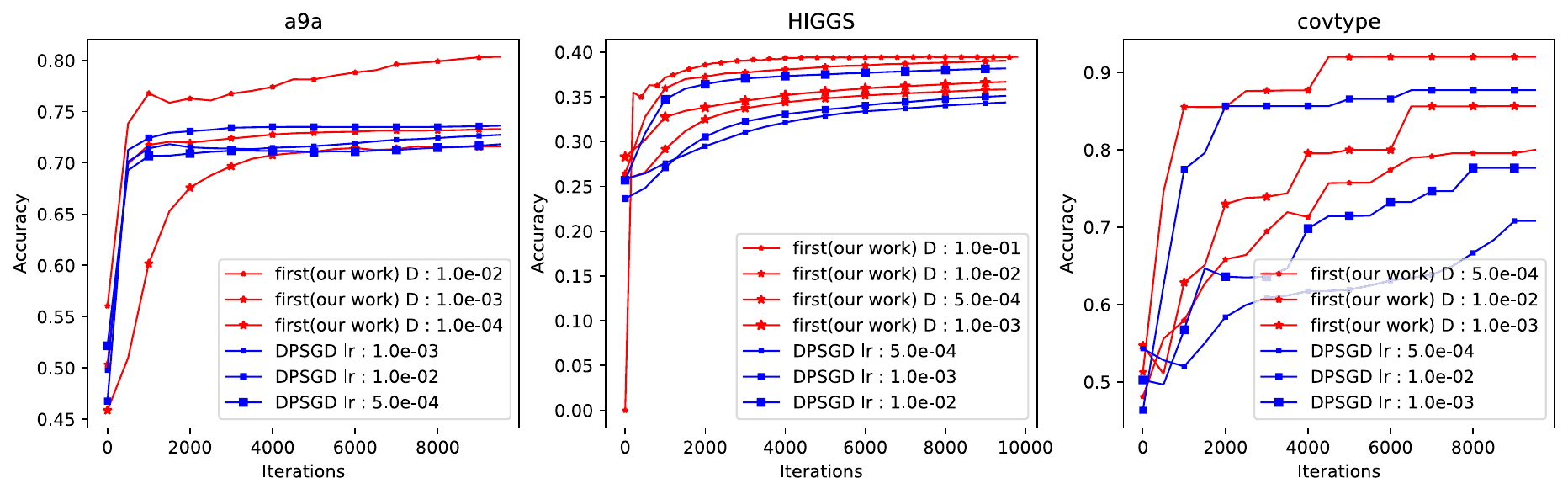}}
\caption{Evaluation of the test accuracy of our algorithm and DPSGD in the first-order setting.}
\label{fig_comp_first}
\end{center}
\end{figure}

In this section, we present the plots of test accuracy comparisons (Figs. \ref{fig_comp_zero}-\ref{fig_comp_first}), described in our experiments.

\end{document}